\documentclass[dvips,sts]{imsart}

\usepackage{natbib}
\usepackage{ams}
\usepackage{graphicx}
\usepackage{epsfig}
\usepackage{soul,color}

\arxiv{math.PR/0000000}

\startlocaldefs

\newcommand{\tr}{^{\prime}}

\def\b#1{\mbox{\boldmath $#1$}}    
\def\bl#1{\mbox{\footnotesize \boldmath {$#1$}}} 

\renewcommand{\th}{\theta}

\newcommand{\be}{\beta}
\newcommand{\de}{\delta}
\newcommand{\la}{\lambda}

\newcommand{\ga}{\gamma}

\newcommand{\diag}{{\rm diag}}    
\newcommand{\pa}{\partial}         
\newcommand{\E}{{\rm E}}         

\newtheorem{example}{Example}

\endlocaldefs

\begin{document}
\begin{frontmatter}
\title{An overview of latent Markov models for longitudinal categorical data}
\runtitle{Latent Markov models}

\begin{aug}
\author{Francesco Bartolucci,
\ead[label=e1]{bart@stat.unipg.it}}
\author{Alessio Farcomeni
\ead[label=e2]{alessio.farcomeni@uniroma1.it}}
\and
\author{Fulvia Pennoni
\ead[label=e3]{fulvia.pennoni@unimib.it}}
%

\runauthor{F. Bartolucci, A. Farcomeni and F. Pennoni}

\address{Department of Economics, Finance and
Statistics, University of Perugia, Via A. Pascoli, 20, 06123
Perugia, Italy, \printead{e1},}
\address{Department of Hygene and Public Health,
Sapienza - University of Rome, Piazzale Aldo Moro, 5,
    00185 Roma, Italy, \printead{e2}}
\address{\and Department of Statistics, University of Milano-Bicocca,
Via Bicocca degli Arcimboldi 8, 20126 Milano, Italy, \printead{e3}.}

\end{aug}

\begin{abstract}
We provide a comprehensive overview of latent Markov (LM) models for
the analysis of longitudinal categorical data. The main assumption
behind these models is that the response variables are conditionally
independent given a latent process which follows a first-order
Markov chain. We first illustrate the basic LM model in which the
conditional distribution of each response variable given the
corresponding latent variable and the initial and transition
probabilities of the latent process are unconstrained. For this
model we also illustrate in detail maximum likelihood estimation
through the Expectation-Maximization algorithm, which may be
efficiently implemented by recursions known in the hidden Markov
literature. We then illustrate several constrained versions of the
basic LM model, which make the model more parsimonious and allow us
to include and test hypotheses of interest. These constraints may be
put on the conditional distribution of the response variables given
the latent process (measurement model) or on the distribution of the
latent process (latent model). We also deal with extensions of LM
model for the inclusion of individual covariates and to multilevel
data. Covariates may affect the measurement or the latent model; we
discuss the implications of these two different approaches according
to the context of application. Finally, we outline methods for
obtaining standard errors for the parameter estimates, for selecting
the number of states and for path prediction. Models and related
inference are illustrated by the description of relevant
socio-economic applications available in the literature.
\end{abstract}

\begin{keyword}
\kwd{EM algorithm} \kwd{Forward-Backward recursions} \kwd{Hidden
Markov models} \kwd{Latent class model} \kwd{Measurement errors}
\kwd{Multilevel model}  \kwd{Panel data} \kwd{Rasch model}
\kwd{Unobserved heterogeneity}
\end{keyword}
\end{frontmatter}

\section{Introduction}

In many applications involving longitudinal data, the interest is
often focused on the evolution of a latent characteristic of a group
of individuals over time, which is measured by one or more
occasion-specific response variables. This characteristic may
correspond, for instance, to the quality-of-life of subjects
suffering from a certain disease, which is indirectly assessed on
the basis of responses to a set of suitably designed items that are
repeatedly administered during a certain period of time.

In the statistical and econometric literatures, several approaches
have been introduced to address the above issue. Among these
approaches, one of the most interesting is  based on  the same
formulation of the hidden Markov model for time series
\citep{baum:petr:66,macd:zucc:97}. The main assumption behind this
approach is that the response variables are conditionally
independent given a latent Markov chain with a finite number of
states. The basic idea behind this assumption, that we will refer to
as the assumption of local independence, is that the latent process
fully explains the observable behavior of a subject; moreover, the
latent state to which a subject belongs at a certain occasion only
depends on the latent state at the previous occasion.

The starting point of the paper is the latent Markov (LM) model.
This model dates back to Wiggins' Ph.D. thesis, \cite{wigg:55}, who
introduced a version of this model based on a homogenous Markov
chain, a single outcome at each occasion, and did not account for
individual covariates. \cite{wigg:55} formulated the model so that a
manifest transition is a mixture of a true change and a spurious
change due to measurement errors in the observed states. See also
\cite{wigg:73} for a deep illustration of this model and some simple
generalizations. It is also worth mentioning \cite{van:del:86},
\cite{vand:lang:90}, \cite{coll:wuga:92}, and \cite{lang:vand:94},
among the first papers dealing with this model.

In the following, we refer to the LM model based on a first-order
Markov chain, non-homogeneous transition probabilities, and no
covariates as the {\em basic LM model}. This model may be used for
univariate or multivariate data; in the second case we observe more
response variables at each occasion. For the basic LM model we
discuss in detail maximum likelihood estimation through the
Expectation-Maximization (EM) algorithm
\citep{baum:et:al:70,demp:lair:rubi:77}, even though we acknowledge
that other estimation methods are available
\citep{arch:titt:02,capp:moul:rydn:05,Kuns:05,turn:08}. For the
implementation of the EM algorithm we illustrate suitable recursions
which allow us to strongly reduce the computational effort.

The paper also focuses on several constrained versions and
extensions of the basic LM model. Constraints have the aim of making
the model more parsimonious and easier to interpret and correspond
to certain hypotheses that may be interesting to test. These
constraints may be posed on the {\em measurement model}, i.e. the
conditional distribution of the response variables given the latent
process, or on the {\em latent model}, i.e. the distribution of the
latent process. About the measurement model, we discuss in detail
\cite{rasch:61} type parameterizations which make the latent states
interpretable in terms of ability or propensity levels. About the
latent model, we outline several simplifications of the transition
matrix, mostly based on constraints of equality between certain
elements of this matrix and/or on the constraint that some elements
are equal to 0. One of the main problems is how to test for these
restrictions. For this aim, we make use of the likelihood ratio (LR)
statistic. It is important to note that, when constrains concern the
transition matrix, the null asymptotic distribution does not
necessarily have an asymptotic chi-squared distribution, but a
distribution of chi-bar-squared type \citep{bart:06}.

The most natural extension of the basic LM model is for the
inclusion of individual covariates. In particular, we describe two
different ways of including such covariates: ({\em i}) in the
measurement model, so that they affect the conditional distribution
of the response variables given the latent process
\citep{bart:farc:09}; ({\em ii}) in the latent model, so that they
affect the initial and the transition probabilities of the Markov
chain \citep{verm:lang:bock:99,bart:penn:fran:07}. Further, we
discuss methods to relax the assumption of local independence and
show multilevel LM models which are suitable when subjects are
collected in clusters. In this context, the model may be formulated
by including fixed parameters to represent the effect of the
clusters on the distribution of the latent process corresponding to
every subject belonging to these clusters; this approach was
followed by \cite{bart:lupp:mont:09}. Otherwise, the model may be
formulated by including random parameters having a discrete
distribution, following in this way an approach similar to that
behind the latent class model, as in \cite{bart:penn:vitt:10}. This
extension is related to the mixed LM model \citep{vand:lang:90} and
to the LM model with random effects \citep{Altman:07}.

Finally, we revise methods for obtaining standard errors for the
model parameters and for selecting the number of latent states. We
also discuss the problem of path prediction through  the Viterbi
algorithm \citep{vite:67,juan:rabi:91}.

An important point to clarify is that, throughout the paper, we
consider the case of categorical response variables because this is
the typical case of application of the LM model. However, it is
straightforward to modify the framework to deal with continuous
outcomes. Most of the theory and estimation methods do not change
substantially.

The paper is organized as follows. In the Section \ref{sec:lm} we
outline the basic LM model and discuss maximum likelihood estimation
for this model based on the EM algorithm. In Section
\ref{sec:modelling} we outline constrained versions of the LM model
based on parsimonious and interpretable parameterizations. In
Section \ref{sec:cov} we illustrate how to deal with individual
covariates, whereas the multilevel extension is presented in Section
\ref{sec:multilevel}. Section \ref{sec:se} deals with standard
errors, selection of the number of states, and path prediction.
Section \ref{applications} illustrates different types of LM model
through various examples involving longitudinal categorical data,
summarizing the results from other papers. The paper ends with a
section where we draw main conclusions and discuss further
developments of the present framework.
\section{Basic latent Markov model and its multivariate version}
\label{sec:lm}
In the following, we illustrate the basic LM model for univariate
categorical data without covariates and in which the latent Markov
chain is of first-order and non-homogenous. We also describe the EM
algorithm for maximum likelihood estimation of this model.
\subsection{Basic formulation of univariate responses}
\label{sec:uni}
Let $\b Y_i = (Y^{(1)}_{i},\ldots,Y_{i}^{(T)})$, $i=1,\ldots,n$, be
a sequence of $T$ categorical response variables with $l$ levels or
categories, coded from $0$ to $l-1$, independently observed over $n$
subjects. Typically, these variables have the same nature, as they
correspond to repeated measurements on the same subjects at
different occasions. However, the approach may be applied to the
case of response variables having a different nature and, possibly,
a different number of categories.

\begin{figure}[ht]\centering
\includegraphics[height=4cm]{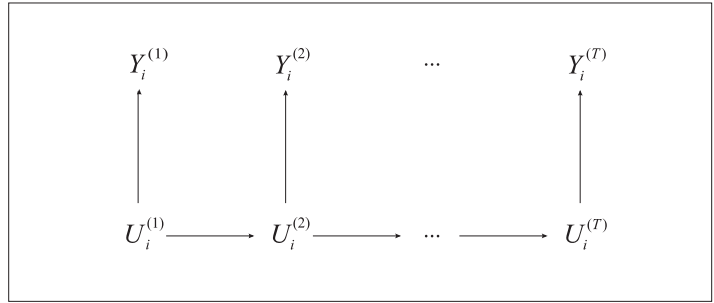}
\caption{{\em Path diagram of the basic LM model for univariate
data.}} \label{fig:basic_LM}
\end{figure}

The main assumption of the basic LM model is that of {\em local
independence}, i.e. for every subject $i$ the response variables in
$\b Y_i$ are conditionally independent given a latent process $\b
U_i=(U^{(1)}_{i},\ldots,U^{(T)}_{i})$. This latent process is
assumed to follow a first-order Markov chain with state space
$\{1,\ldots,k\}$. Then, for all $t>2$, the latent variable
$U^{(t)}_{i}$ is conditionally independent of
$U^{(1)}_{i},\ldots,U^{(t-2)}_{i}$ given $U^{(t-1)}_{i}$. See Figure
\ref{fig:basic_LM} for an illustration via path diagram.

The latent process can be used to model a real unidimensional latent
trait evolving over time, or just as a parsimonious device to allow
for a time non-homogeneous distribution for the sequence $\b Y_i$.
The Markov assumption is seldom found to be restrictive, and is
easily interpretable. A further interpretation of the model is given
by the {\it measurement error} framework: the unobservable outcome
$U^{(t)}_{i}$ is observed with measurement error as $Y^{(t)}_{i}$.
Therefore, it may be seen as an extension of a standard Markov model
\citep{ande:54}; for a discussion on this point see \cite[Chapter
4]{wigg:73}.

Parameters of the model are the conditional response probabilities
$\phi_{y|u}^{(t)}=p(Y^{(t)}_{i}=y|U^{(t)}_{i}=u)$, $t=1,\ldots,T$,
$u=1,\ldots,k$, $y=0,\ldots,l-1$, the initial probabilities $\pi_{u}
= p(U^{(1)}_{i}=u)$, $u=1,\ldots,k$, and the transition
probabilities $\pi_{v|u}^{(t)} = p(U^{(t)}_{i}=v|U^{(t-1)}_{i}=u)$,
$t=2,\ldots,T$, $u,v=1,\ldots,k$. Note that all these probabilities
do not depend on $i$ since, in its basic version, the model does not
account for individual covariates.

On the basis of the above parameters, the distribution of $\b U_i$
may be expressed as
\begin{equation}
p(\b U_i=\b
u)=\pi_{u^{(1)}}\prod_{t>1}\pi_{u^{(t)}|u^{(t-1)}}^{(t)},
\label{eq:latprob}
\end{equation}
where $\b u = (u^{(1)},\ldots,u^{(T)})$. Moreover, the conditional
distribution of $\b Y_i$ given $\b U_i$ may be expressed as
\[
p(\b Y_i=\b y|\b U_i=\b u) = \prod_t \phi_{y^{(t)}|u^{(t)}}^{(t)},
\]
and, consequently, for the {\em manifest distribution} of $\b Y_i$
we have
\begin{eqnarray}
\label{eq:manifest_probability}\hspace*{0.75cm}f(\b y) &=& p(\b Y_i=\b y)=\sum_{\bl u}f(\b u,\b y)=\\
&=&\sum_{u^{(1)}}\phi_{y^{(1)}|u^{(1)}}^{(1)}\pi_{u^{(1)}}
\sum_{u^{(2)}}\phi_{y^{(2)}|u^{(2)}}^{(2)}\pi_{u^{(2)}|u^{(2)}}^{(2)}\cdots
\sum_{u^{(T)}}\phi_{y^{(T)}|u^{(T)}}^{(T)}\pi_{u^{(T)}|u^{(T-1)}}^{(T)},\nonumber
\end{eqnarray}
where $\b y=(y^{(1)},\ldots,y^{(T)})$ and $f(\b u,\b y)=p(\b Y_i=\b
y|\b U_i=\b u)p(\b U_i=\b u)$. It is important to note that
computing $f(\b y)$ as expressed above involves a sum extended to
all the possible $k^T$ configurations of the vector $\b u$;  this
typically requires a considerable computational effort.

In order to efficiently compute the probability $f(\b y)$, we can use a
forward recursion \citep{baum:et:al:70} for obtaining $q^{(t)}(u,\b
y)=p(U^{(t)}_{i}=u,Y^{(1)}_{i}=y^{(1)},\ldots,Y^{(t)}_{i}=y^{(t)})$
for $t=1,\ldots,T$. We then have
\[
f(\b y)=\sum_u q^{(T)}(u,\b y).
\]
In particular, given $q^{(t-1)}(u,\b y)$, $u=1,\ldots,k$, the $t$-th
iteration of the recursion, $t=2,\ldots,T$, consists of computing
\begin{equation}
q^{(t)}(v,\b y) = \sum_u q^{(t-1)}(u,\b
y)\pi^{(t)}_{u|v}\phi^{(t)}_{y^{(t)}|u^{(t)}},\quad
v=1,\ldots,k,\label{eq:rec}
\end{equation}
starting with $q^{(1)}(u,\b y)=\pi_u\phi^{(1)}_{y^{(1)}|u}$.

The above recursion may be simply implemented by using the matrix
notation \citep{bart:penn:fran:07}. Let $\b q^{(t)}(\b y)$ denote
the column vector with elements $q^{(t)}(u,\b y)$, $u=1,\ldots,k$,
and in similar way define the initial probability vector $\b\pi$
with elements $\pi_u$ and the conditional probability vector
$\b\phi^{(t)}_y$ with elements $\phi^{(t)}_{y|u}$. Also let
$\b\Pi^{(t)}$ denote the transition probability matrix with elements
$\pi^{(t)}_{v|u}$, $u,v=1,\ldots,k$, with $u$ running by row and $v$
by column. We then have
\begin{equation}
\b q^{(t)}(\b y) = \left\{\begin{array}{ll}
\diag[\b\phi^{(1)}_{y^{(1)}}]\b\pi, & \mbox{ if } t=1,\\
\diag[\b\phi^{(t)}_{y^{(t)}}][\b\Pi^{(t)}]\tr\b q^{(t-1)}(\b y), & \mbox{ otherwise},\\
\end{array}\right.\label{eq:forward_recursion}
\end{equation}
and at the end $f(\b y) =\b q^{(T)}(\b y)\tr\b 1$, where $\b 1$
denotes a column vector of ones of suitable dimension. In
implementing this recursion, attention must be payed to the case of
large values of $T$, because during the recursion the probabilities
$q^{(t)}(u,\b y)$ could become very small; see \cite{Scott:02}.
\subsection{Multivariate version}
\label{multiLMmodels}
In the multivariate case we observe a vector $\b Y^{(t)}_{i} =
(Y^{(t)}_{i1},\ldots,Y^{(t)}_{ir})$ of $r$ response variables for
every subject $i$ and occasion $t$. For $j=1,\ldots,r$, the response
variable $Y^{(t)}_{ij}$ has $l_j$ levels coded from $0$ to $l_j-1$.
We also denote by $\b Y_i$ the response vector made of the union of
the vectors $\b Y^{(t)}_{i}$ for $t=1,\ldots,T$.

\begin{figure}[ht]\centering
\includegraphics[height=4cm]{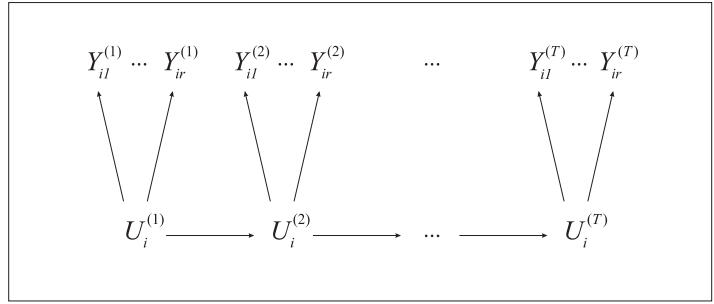}
\caption{\em Path diagram of the basic LM model for multivariate
data.} \label{fig:multivariate_LM}
\end{figure}

The LM model presented in Section \ref{sec:uni} may be naturally
formulated for multivariate data by considering an extended version
of the assumption of local independence. Further to assuming that
for every $i$ the vectors $\b Y^{(t)}_{i}$, $t=1,\ldots,T$, are
conditionally independent given $\b U_i$, we assume that the
response variables in each vector $\b Y^{(t)}_{i}$ are conditionally
independent given $U^{(t)}_{i}$. The resulting model is represented
by the path diagram in Figure \ref{fig:multivariate_LM}

The model assumptions imply that
\begin{equation}
p(\b Y_i=\b y|\b U_i=\b u)=\prod_t p(\b Y^{(t)}_{i}=\b
y^t|U^{(t)}_{i}=u^{(t)}),\label{eq:prob_cond_multi}
\end{equation}
where $\b y$ is made of the subvectors $\b
y^{(t)}=(y^{(t)}_{1},\ldots,y^{(t)}_{r})$ and
\begin{equation}
p(\b Y^{(t)}_{i}=\b
y^{(t)}|U^{(t)}_{i}=u^{(t)})=\prod_j\phi^{(t)}_{j,y^{(t)}_{j}|u^{(t)}},\label{eq:joint_prob}
\end{equation}
with $\phi^{(t)}_{j,y|u}=p(Y^{(t)}_{ij}=y|U^{(t)}_{i}=u)$. The
manifest probability $f(\b y)$ has then the same expression as in
(\ref{eq:manifest_probability}), with $p(\b Y_i=\b y|\b U_i=\b u)$
computed as in (\ref{eq:prob_cond_multi}). This manifest probability
may be computed by exploiting the same recursion as in
(\ref{eq:forward_recursion}), with $\b\phi^{(t)}_{y^{(t)}}$
substituted by a vector of the same dimension with elements $p(\b
Y^{(t)}_{i}=\b y^{(t)}|U^{(t)}_{i}=u^{(t)})$, $u=1,\ldots,k$.
\subsection{Maximum likelihood estimation}
 \label{sec:EM}
For an observed sample of $n$ subjects, let $\b y_i$ denote the
(univariate or multivariate) response configuration provided by
subject $i$. The log-likelihood of the LM model may be expressed as
\begin{equation}
\ell(\b\th) = \sum_i \log[f(\b y_i)],\label{eq:lk}
\end{equation}
where $\b\th$ is the vector of all model parameters arranged in a
suitable way. An equivalent form is
\begin{equation}
\ell(\b\th) = \sum_{\bl y} n_{\bl y}\log[f(\b y)],\label{eq:lk_freq}
\end{equation}
where $n_{\bl y}$ denotes the frequency of the response
configuration $\b y$ in the sample. In absence of covariates, using
(\ref{eq:lk_freq}) for computing $\ell(\b\th)$ is more efficient
since the sum $\sum_{\bl y}$ may be restricted to all response
configurations $\b y$ observed at least once; we then adopt this
formulation. We estimate $\b\th$ by maximizing the log-likelihood
$\ell(\b\th)$. This may be easily done by the
Expectation-Maximization (EM) algorithm
\citep{baum:et:al:70,demp:lair:rubi:77}.

The EM algorithm is based on the concept of {\em complete data},
which in our case consist of the response configuration $\b y_i$
(incomplete data) and the configuration of the latent process $\b
u_i$ (missing data) for every subject $i$. An equivalent way to
represent the complete data, which is more coherent with
(\ref{eq:lk_freq}), is by the frequencies $m_{\bl u\bl y}$ of the
contingency table in which the subjects are cross-classified
according to the latent configuration $\b u$ and the response
configuration $\b y$. Under this formulation, the complete data
likelihood has logarithm
\[
{\ell^*}(\b\th)=\sum_{\bl u}\sum_{\bl y}m_{\bl u\bl y}\log[f(\b u,\b
y)].\label{eq:complk}
\]
After some simple algebra, in the univariate case we have

\begin{eqnarray}
{\ell^*}(\b\th)&=&\sum_t\sum_u\sum_y {\tilde{a}}^{(t)}_{uy}\log[\phi^{(t)}_{y|u}]+\label{eq:comp_lik_simpl}\\
&+&\sum_u a^{(1)}_{u}\log(\pi_u)+\sum_{t>1}\sum_u\sum_v
a^{(t)}_{uv}\log[\pi_{v|u}^{(t)}],\nonumber
\end{eqnarray}
where $a^{(t)}_{u}$ is the number of subjects that at occasion $t$
are in latent state $u$; with reference to the same occasion $t$,
$a^{(t)}_{uv}$ is the number of subjects that move from latent state
$u$ to latent state $v$, and ${\tilde{a}}^{(t)}_{uy}$ is the number
of subjects that are in latent state $u$ and respond by $y$. Note
that ${\ell^*}(\b\th)$ is made of three components that may be
separately maximized; see \cite{bart:06} for details. Also note that
a simplification similar to (\ref{eq:comp_lik_simpl}) holds for the
multivariate case.

The frequencies $a^{(t)}_{u}$, $a^{(t)}_{uv}$ and
$\tilde{a}^{(t)}_{uy}$ above are obviously unknown. Then, the EM
algorithm proceeds by alternating the following two steps until
convergence in $\ell(\b\th)$:
\begin{itemize}
\item{\bf E-step}: it consists of computing the expected value of
each unknown frequency in (\ref{eq:comp_lik_simpl}) given the
observed data and the current value of the parameters, so as to
obtain the expected value of ${\ell^*}(\b\th)$. The expected values
of these frequencies are obtained as:
\begin{eqnarray*}
\hat{a}^{(t)}_{u}&=&\sum_{\bl y}n_{\bl y}r^{(1)}(u|\b y),\\
\hat{a}^{(t)}_{uv}&=&\sum_{\bl y}n_{\bl y}\sum_{t>1}r^{(t)}(u,v|\b y),\\
\hat{\tilde{a}}^{(t)}_{uy}&=&\sum_{\bl y}n_{\bl
y}I(y^{(t)}=y)r^{(t)}(u|\b y),
\end{eqnarray*}
where $I(\cdot)$ is the indicator function, $r^{(t)}(u|\b
y)=p(U^{(t)}_{i}=u|\b Y_i=\b y)$, and $r^{(t)}(u,v|\b
y)=p(U^{(t-1)}_{i}=u,U^{(t)}_{i}=v|\b Y_i=\b y)$.
\item{\bf M-step}: it consists of updating the estimate of $\b\th$
by maximizing the expected value of ${\ell^*}(\b\th)$ obtained as
above. Explicit solutions are available at this aim. In particular,
we have:
\begin{itemize}
\item{\em Conditional probabilities of response}:
${\displaystyle\phi_{y|u}^{(t)}=\frac{\hat{\tilde{a}}^{(t)}_{uy}}{\hat{a}^{(t)}_u}}$,
$t=1,\ldots,T$, $u=1,\ldots,k$, $y=0,\ldots,l-1$.
\item{\em Initial probabilities}: ${\displaystyle\pi_u = \frac{\hat{a}^{(1)}_u}{n}}$, $u=1,\ldots,k$.
\item{\em Transition probabilities}:
${\displaystyle\pi_{v|u}^{(t)} =
\frac{\hat{a}^{(t)}_{uv}}{\hat{a}^{(t-1)}_u}}$, $t=2,\ldots,T$,
$u,v=1,\ldots,k$.
\end{itemize}
\end{itemize}
A rule similar to the above one must be applied to update the
parameters $\phi_{j,y|u}^{(t)}$ in the multivariate case.

A final point is how to compute in an efficient way the posterior
probabilities $r^{(t)}(u,v|\b y)$ on the basis of which we obtain
the posterior probabilities $r^{(t)}(u|\b y)$ via marginalization.
First of all consider the probabilities $q^{*(t)}(u,\b
y)=p(Y^{(t+1)}_{i}=y^{(t+1)},\ldots,Y^{(T)}_{i}=y^{(T)}|U^{(t)}=u)$.
For $t=1,\ldots,T-1$, these probabilities may be computed by the
backward recursion
\[
q^{*(t)}(u,\b y) =\sum_v q^{*(t+1)}(v,\b
y)\pi^{(t+1)}_{v|u}\phi^{(t+1)}_{y_{t+1}|v},\quad u=1,\ldots,k,
\]
initialized with  $q^{*(T)}(v,\b y)=1$, $v=1,\ldots,k$ \citep[Sec.
2.2]{baum:et:al:70,levi:rabi:sond:83,macd:zucc:97}. Then, for
$t=2,\ldots,T$, we obtain the posterior probabilities above as
\[
r^{(t)}(u,v|\b y) = \frac{q^{(t-1)}(u,\b
y)\pi^{(t)}_{v|u}\phi^{(t)}_{y^{(t)}|u}q^{*(t)}(v,\b y)}{f(\b
y)},\quad u,v=1,\ldots,k.
\]
See \cite{bart:besa:09} for an alternative recursion to compute posterior
probabilities.

Even in this case, the recursion can be efficiently implemented by
the matrix notation \citep{bart:penn:fran:07}. Let $\b q^{*(t)}(\b
y)$ be the vector with elements $q^{*(t)}(u,\b y)$, $u=1,\ldots,k$.
This vector may computed as
\begin{equation}
\b q^{*(t)}(\b y) = \left\{\begin{array}{ll}
\b 1, & \mbox{ if } t=T,\\
\b\Pi^{(t+1)}\diag(\b\phi_{y^{(t+1)}}^{(t+1)})\b q^{*(t+1)}(\b y), &
\mbox{ otherwise.}
\end{array}\right.\label{eq:rec2}
\end{equation}
Then, the $k\times k$ matrix $\b R^{(t)}(\b y)$, with elements
$r^{(t)}(u,v)$ arranged by letting $u$ run by row and $v$ by column,
is obtained as
\begin{equation}
\b R^{(t)}(\b y) = \frac{\diag[\b q^{(t-1)}(\b
y)]\b\Pi^{(t)}\diag[\b\phi^{(t)}_{y^{(t)}}]\diag[\b q^{*(t)}(\b
y)]}{f(\b y)},\quad t=2,\ldots,T.\label{eq:defR}
\end{equation}
A similar implementation holds for the multivariate case by
substituting each vector $\b\phi^{(t)}_{y^{(t)}}$ with the vector
having elements $p(\b Y^{(t)}_{i}=\b y^{(t)}|U^{(t)}_{i}=u)$.

As typically happens for latent variable and mixture models, the
likelihood function may be multimodal. In particular, the EM
algorithm could converge to a mode of the likelihood which does not
correspond to the global maximum. In order to increase the chance of
reaching the global maximum, the EM needs to be initialized in a
proper way. Further, the value at convergence may be compared with
values obtained starting from initial values randomly chosen; for a
similar multi-start strategy for mixture models see
\cite{berchtold:04}.  Inference is then based on the solution
corresponding to the highest value of the likelihood at convergence.
We take this solution as the maximum likelihood estimates of $\b
\th$ donated by $\hat{\b \th}$.

\section{Constrained versions of the LM model}
\label{sec:modelling}
In this section, we show how we can adopt a more parsimonious
parameterization on the basis of constraints which typically
correspond to hypotheses of interest. These constraints may concern:
\begin{itemize}
\item{\em measurement model}: the distribution of the response
variables given the latent process which depends on the conditional
response probabilities $\phi_{y|u}^{(t)}$ or  $\phi_{j,y|u}^{(t)}$;
\item{\em latent model}: the distribution of the latent process
which depends on the initial probabilities $\pi_u$ and the
transition probabilities $\pi_{v|u}^{(t)}$.
\end{itemize}

We also discuss maximum likelihood estimation and the likelihood
ratio testing of constraints on the model parameters.
\subsection{Constraints on the measurement model}\label{sec:modelling1}
\subsubsection{Univariate case.}
With only one response variable at each occasion, a sensible
constraint is that
\begin{equation}
\phi_{y|u}^{(t)}=\phi_{y|u},\quad
t=1,\ldots,T,\:u=1,\ldots,k,\:y=0,\ldots,l-1.\label{eq:same_phi}
\end{equation}
This constraint corresponds to the hypothesis that the distribution
of the responses only depends on the corresponding latent variable
and so there is no dependence of this distribution on time.

More sophisticated constraints may be formulated in the form
\begin{equation}
\b\eta_u^{(t)}=\b Z_u^{(t)}\b\be,\label{eq:glm}
\end{equation}
where $\b\eta_u^{(t)}=\b
g(\phi^{(t)}_{0|u},\ldots,\phi^{(t)}_{l-1|u})$, with $\b g(\cdot)$
being a suitable link function, and $\b Z_u^{(t)}$ a design matrix.
With binary response variables, a natural link function is of logit
type, so that
\[
\eta_u^{(t)}=\log\frac{\phi^{(t)}_{1|u}}{\phi^{(t)}_{0|u}},
\]
or a probit link function based on the normal distribution function.
With response variables having more than two categories, a natural
choice is that of a multinomial logit link function, so that
$\b\eta_u^{(t)}$ has elements
\[
\log\frac{\phi^{(t)}_{y|u}}{\phi^{(t)}_{0|u}},\quad y=1,\ldots,l-1.
\]
With ordinal response variables, more sensible logits are of global
or continuation type; in the first case $\b\eta_u^{(t)}$ has
elements
\[
\log\frac{p(Y^{(t)}_{i}\geq y|U^{(t)}_{i}=u)}{p(Y^{(t)}_{i}<
y|U^{(t)}_{i}=u)}=\log\frac{\phi^{(t)}_{y|u}+\cdots+\phi^{(t)}_{l-1|u}}
{\phi^{(t)}_{0|u}+\cdots+\phi^{(t)}_{y-1|u}},\quad y=1,\ldots,l-1,
\]
and in the second it has elements
\[
\log\frac{p(Y^{(t)}_{i}\geq
y|U^{(t)}_{i}=u)}{p(Y^{(t)}_{i}=(y-1)|U^{(t)}_{i}=u)}=
\log\frac{\phi^{(t)}_{y|u}+\cdots+\phi^{(t)}_{l-1|u}}
{\phi^{(t)}_{0|u}},\quad y=1,\ldots,l-1;
\]
see \cite{colo:forc:01} and the references therein for a
comprehensive review of these types of logit. Obviously, constraint
(\ref{eq:same_phi}) may be included into (\ref{eq:glm}) by requiring
$\b Z_u^{(t)}=\b Z_u$, $t=1,\ldots,T$, $u=1,\ldots,k$. However,
using constraint (\ref{eq:same_phi}) per se leads to a model which
is easier to estimate and does not require to rely on a link
function.

An interesting example about parameterization (\ref{eq:glm}) is
given in the following example.\vskip5mm

\begin{example}{\bf - LM Rasch model.} This is a version of the
\cite{rasch:61}  model for binary data, based on the assumption that
the ability evolves over time, which is formulated by assuming
\begin{equation}
\eta_{t|u} =
\log\frac{\phi^{(t)}_{1|u}}{\phi^{(t)}_{0|u}}=\xi_u-\psi^{(t)},\qquad
t=1,\ldots,s,\quad u=1,\ldots,c.\label{eq:rasch}
\end{equation}
This formulation makes sense for data derived from the
administration of a set of $T$ test items to a group of $n$
subjects, a situation that frequently arises in psychological and
educational measurement. In this case, $\xi_u$ may be interpreted as
the ability level of subjects in latent state $u$, whereas
$\psi^{(t)}$ may be interpreted as the difficulty level of item $t$.
Note that in this case we may have response variables of different
nature, since they correspond to different items. So, strictly
speaking, we are not in a longitudinal context in which the same
response variable is repeatedly observed. Nevertheless, the model
makes sense as an alternative to the latent class Rasch model
\citep{deLee:ver:86,lind:clo:gre:91} and to test violation of the
latter. For a detailed description see \cite{bart:penn:lupp:08}.
Finally note that parametrizations similar to (\ref{eq:rasch}) may
be adopted with ordinal variables on the basis of logits of global
or continuation types.
\end{example}

With binary responses, an advantage of  parameterization
(\ref{eq:rasch}) is that it implies the hypothesis of {\em
monotonicity}, i.e. the latent states can always be ordered so that
the probability of success increases with the label of the latent
class, in symbols
\begin{equation}
\phi^{(t)}_{1|1}\leq\cdots\leq\phi^{(t)}_{1|k},\quad
t=1,\ldots,T.\label{eq:monotonocity}
\end{equation}
Note that this constraint may be included in the model regardless of
a specific parameterization as the one above. In particular,
\cite{bart:06}  allows for inequality constraints in the general
form $\b K\b\be\geq\b 0$, where $\b 0$ is a vector of zeros of
suitable dimension, so that we can obtain ordered latent states in
non-parametric way; see \cite{bart:for:05} and the references
therein.
\subsubsection{Multivariate case}
In the multivariate case, the parameterizations suggested above may
be adopted for each single response variable; then, the conditional
distribution of each vector $\b Y^{(t)}_{i}$ given $U^{(t)}_{i}$ is
still obtained as in (\ref{eq:joint_prob}).

In particular, constraint (\ref{eq:same_phi}) becomes
\begin{equation}
\hspace*{0.75cm} \phi^{(t)}_{j,y|u}=\phi_{j,y|u},\quad
j=1,\ldots,r,\:t=1,\ldots,T,\:u=1,\ldots,k,\:y=0,\ldots,l_j-1.\label{eq:same_phi_multi}
\end{equation}
Moreover, we can adopt a specific link function to parameterize the
conditional distribution of each response variable as in
(\ref{eq:glm}). We let $\b\eta_{j,u}^{(t)} = \b
g_h(\phi_{j,0|u}^{(t)},\ldots,\phi_{j,l_j-1|u}^{(t)})$, so that we
formulate the model as
\begin{equation}
\b\eta_{j,u}^{(t)}=\b Z_{j,u}^{(t)}\b\be.\label{eq:glm_multi}
\end{equation}
\subsection{Constraints on the latent model}
\label{sec:modelling2}
In absence of individual covariates, no many interesting constraints
may be expressed on the initial probabilities $\pi_u$. The only
constraint that may be of interest is that of uniform initial
probabilities
\begin{equation}
\pi_u=1/k,\quad u=1,\ldots,k,\label{eq:uniform_pi}
\end{equation}
meaning that at the beginning of the survey each latent state has
the same proportion of subjects.

More interesting constraints may be expressed on the transition
probabilities. These constraints allow us to strongly reduce the
number of parameters of the model. A simple constraint to express is
that the Markov chain is time homogenous, that is
\begin{equation}
\pi_{v|u}^{(t)}=\pi_{v|u},\quad t=2,\ldots,T,\quad
u,v=1,\ldots,k.\label{eq:time_homogenous}
\end{equation}
Note that we can also adopt a constraint of partial homogeneity; for
instance, \cite{bart:penn:fran:07}  adopted two different transition
matrices, one until occasion $\bar{t}$ and the other for transitions
after this occasion, that is
\begin{equation}
\pi_{v|u}^{(t)}=\left\{\begin{array}{ll}
\bar{\pi}_{v|u}^{(1)}, & t=2,\ldots,\bar{t},\\
\bar{\pi}_{v|u}^{(2)}, & t=\bar{t}+1,\ldots,T,
\end{array}\right.\label{eq:partial_homo}
\end{equation}
with $\bar{t}$ between 2 and $T$ and $u,v=1,\ldots,k$.

More sophisticated parametrizations may be formulated by a linear or
a generalized linear model on the transition probabilities. Linear
models have the advantage of permitting to express the constraint
that certain probabilities are equal to 0, so that transition
between two given states is not possible. In particular,
\cite{bart:06} considered a formulation that in our case may be
expressed as
\begin{equation}
\b\rho_u^{(t)}=\b W_u^{(t)}\b\de,\label{eq:linear_model}
\end{equation}
where $\b\rho_u^{(t)}$ is the column vector containing the elements
of the $u$-th row of the $t$-th transition matrix, apart from the
diagonal element, i.e. $\pi_{v|u}^{(t)}$, $v=1,\ldots,k$, $v\neq u$,
and $\b W_u^{(t)}$ is a suitable design matrix. In order to ensure
that all the transition probabilities are non-negative, we have to
impose suitable restrictions on the parameter vector $\b\de$. Due to
these restrictions, estimation may be more difficult and we are not
in a standard inferential problem; see Section
\ref{sec:testing_constraints}.

\begin{example}
{\bf - Transition matrices that may be formulated by the linear
model.} The simplest constraint is that all the off-diagonal
elements of the transition matrix $\b\Pi^{(t)}$ are equal to each
other; with $k=3$, for instance, if we also assume homogeneity, we
have
\begin{equation}
\b\Pi^{(t)}=\pmatrix{1-2\de^{(t)} & \de^{(t)} & \de^{(t)}\cr
\de^{(t)} & 1-2\de^{(t)} & \de^{(t)}\cr \de^{(t)} & \de^{(t)} &
1-2\de^{(t)}},\quad t=2,\ldots,T.\label{eq:trans_only_one}
\end{equation}
This constraint may be formulated by letting all design matrices $\b
W_u^{(t)}$ to be simply equal to $\b 1$. A less stringent constraint
is that each transition matrix is symmetric, so that the probability
of transition from latent state $u$ to latent state $v$ is the same
as that of the reverse transition:
\[
\b\Pi^{(t)}= \pmatrix{1-(\de_1^{(t)}+\de_2^{(t)}) & \de_1^{(t)} &
\de_2^{(t)}\cr \de_1^{(t)} & 1-(\de_1^{(t)}+\de_3^{(t)}) &
\de_3^{(t)}\cr \de_2^{(t)} & \de_3^{(t)} &
1-(\de_2^{(t)}+\de_3^{(t)})}.
\]
Finally, when the latent states are ordered in a meaningful way by
assuming, for instance, that (\ref{eq:monotonocity}) holds, it may
be interesting to formulate the hypothesis that a subject in latent
state $u$ may move only to latent state $v=u+1,\ldots,k$. With
$k=3$, for instance, we have
\[
\b\Pi^{(t)}=\pmatrix{1-(\de_1^{(t)}+\de^{(t)}_2) & \de^{(t)}_1 &
\de^{(t)}_2\cr 0 & 1-\de^{(t)}_3 & \de^{(t)}_3\cr 0 & 0 & 1}.
\]
A first  example  of such restriction was provided by
\cite{coll:wuga:92}, in which latent states represent ordered
developmental states. According to the underlying developmental
psychology theory, children may make a transition to a next stage
but will never return to a previous stage.
\end{example}

An alternative parameterization is based on using a suitable link
function for each row of the transition matrix. The model may be
then formulated as
\begin{equation}
\b\la_u^{(t)}=\b V_u^{(t)}\b\de,\label{eq:glm_transition}
\end{equation}
with $\b\la^{(t)}_u=\b g_u(\pi_{1|u}^{(t)},\ldots,\pi_{k|u}^{(t)})$
and $\b V_u^{(t)}$ being a suitable design matrix. A possible link
function is based on logits with respect to the diagonal element, so
that $\b\la^{(t)}_u$ has $k-1$ equal to
\begin{equation}
\log\frac{\pi_{v|u}^{(t)}}{\pi_{u|u}^{(t)}},\quad
v=1,\ldots,k,\:v\neq u. \label{eq:logit_transition}
\end{equation}
An alternative parameterization, which makes sense with ordered
latent states, is based on global logits, so that the elements of
each vector $\b\la^{(t)}_u$ are
\begin{equation}
\log\frac{\pi_{v|u}^{(t)}+\cdots+\pi_{k|u}^{(t)}}
{\pi_{1|u}^{(t)}+\cdots+\pi_{v-1|u}^{(t)}},\quad
v=2,\ldots,k.\label{eq:global_transition}
\end{equation}
Many other link functions may be formulated, such as the same
adopted within the ordered probit model.

It is worth noting that we can combine a parametrization of type
(\ref{eq:glm_transition}) with the constraint that certain
transition probabilities are equal to 0. In this case the link
function $\b g_u(\cdot)$ must be applied to only those elements in
the $u$-th row of the $t$-th transition matrix that are not
constrained to be equal to 0. In general, the size of each vector
$\b\la_u^{(t)}$ is equal to the number of these elements minus 1 and
the design matrices $\b V_u^{(t)}$ in (\ref{eq:glm_transition}) need
to be defined accordingly. The following example clarifies this
case.

\begin{example}{\bf - Tridiagonal transition matrix.}
A strong reduction of the number of parameters may be achieved by
assuming that each transition matrix is tridiagonal, so that
transition from state $u$ is only allowed to state $v=u-1,u+1$; with
$k=4$, for instance, we have
\begin{equation}
\b\Pi^{(t)}=\pmatrix{\pi_{1|1}^{(t)} & \pi_{2|1}^{(t)} & 0 & 0\cr
\pi_{1|2}^{(t)} & \pi_{2|2}^{(t)} & \pi_{3|2}^{(t)} & 0\cr 0 &
\pi_{2|3}^{(t)} & \pi_{3|3}^{(t)} & \pi_{4|3}^{(t)}\cr 0 & 0 &
\pi_{3|4}^{(t)} & \pi_{4|4}^{(t)}}.\label{eq:triang_transition}
\end{equation}
This constraint makes sense only if latent states are suitably
ordered. In this case, it may combined with a parameterization based
on logits of type (\ref{eq:logit_transition}) by assuming, for
instance,
\[
\log\frac{\pi^{(t)}_{v|u}}{\pi^{(t)}_{u|u}}=\xi_u+\psi_v,\quad
t=2,\ldots,T,\:u=1,\ldots,v,
\]
where $v=2$ for $u=1$, $v=k-1$ for $u=k$ and $v=u-1,u+1$ for
$u=2,\ldots,k-1$. This parameterization may be still expressed as in
(\ref{eq:glm_transition}); note that in this case the vector
$\b\la_u^{(t)}$ containts only one logit for $u=1,k$ and 2 logits
for $u=2,\ldots,k-1$.
\end{example}
\subsection{Maximum likelihood estimation}\label{sec:EM_cov}
Under the constraints illustrated in Sections \ref{sec:modelling1}
and \ref{sec:modelling2},  maximum likelihood estimation of the
parameters of the LM model is carried out by the EM algorithm
illustrated in Section \ref{sec:EM} in which we only have to modify
the M-step according to the constraint of interest. In the
following, we describe in detail the required adjustments. In any
case, along the same lines as \cite{shi:zhen:guo:05}, it is possible
to prove that the observed log-likelihood $\ell(\b\th)$ at each EM
iteration is not decreasing and, therefore, the algorithm converges
to a local maximum of the log-likelihood.

As usual, in order to increase the chance of converging to the
global maximum, we recommend a multi-start strategy as that outlined
at the end of Section \ref{sec:EM}.
\subsubsection{Constraints on the measurement
model.}\label{sec:EM_measurement_cov}
Under constraint (\ref{eq:same_phi}), the parameters of the
conditional distribution of the response variables given the latent
process are updated as
\[
\phi_{y|u}=\frac{\sum_t\hat{\tilde{a}}^{(t)}_{uy}}{\sum_t\hat{a}^{(t)}_u},\quad
u=1,\ldots,k,\: y=0,\ldots,l-1.
\]
A similar rule must be applied in order to update the parameters
$\phi_{j,y|u}$ of this distribution in the multivariate case, when
constraint (\ref{eq:same_phi_multi}) is assumed.

Under a parameterization of type (\ref{eq:glm}), this  conditional
distribution depends on the parameter vector $\b\be$. This parameter
vector is updated at the M-step by maximizing the corresponding
component of the expected value of the complete data log-likelihood,
i.e.
\[
\hat{\ell}^*_1(\b\be)=\sum_t\sum_u\sum_y
\hat{\tilde{a}}^{(t)}_{uy}\log[\phi^{(t)}_{y|u}].
\]
This may be done by standard iterative algorithms,  such as the
Fisher-scoring, which is based on the score vector and the expected
information matrix corresponding to $\hat{\ell}^*_1(\b\be)$. These
quantities have specific expressions depending on the adopted link
function which, in turn, depends on the nature of the response
variables.

Consider, for instance, the case of binary response variables in
which a logit link function is adopted, where the design matrix $\b
Z_u^{(t)}$ specializes into the row vector $[\b z_u^{(t)}]\tr$. The
score vector and the information matrix have, respectively, the
following expressions:
\begin{eqnarray*}
\hat{\b s}_1^*(\b\be)&=&\sum_t\sum_u\sum_y \hat{\tilde{a}}^{(t)}_{uy}(y-\phi^{(t)}_{1|u})\b z_u^{(t)},\\
\hat{\b F}_1^*(\b\be)&=&\sum_t\sum_u
\hat{a}^{(t)}_{u}\phi^{(t)}_{0|u}\phi^{(t)}_{1|u}\b z_u^{(t)}[\b
z_u^{(t)}]\tr.
\end{eqnarray*}
Similar expressions may be obtained for the other cases; see
\cite{bart:06}  for a general discussion which also concerns the use
of constraints of type $\b K\b\be\geq\b 0$. The same iterative
algorithms may be used in the multivariate case, when a
parameterization of type (\ref{eq:glm_multi}) is adopted.
\subsubsection{Constraints on the latent model.}\label{sec:EM_latent_cov}
Under the constraint of uniform initial probabilities, the M-step
simply skips the update of this parameters, whose values are fixed
as in (\ref{eq:uniform_pi}).

About the transition probabilities, when the Markov chain is assumed
to be time homogenous, see equation (\ref{eq:time_homogenous}),
these probabilities are updated as follows:
\[
\pi_{v|u}=\frac{\sum_{t>1}\hat{a}^{(t)}_{uv}}{\sum_{t>1}\hat{a}^{(t-1)}_{u}},
\quad u,v=1,\ldots,k.
\]
A similar rule may be applied in the partial homogeneity case to
update the parameters $\bar{\pi}^{(1)}_{v|u}$ and
$\bar{\pi}^{(2)}_{v|u}$. For the parameters of the first type, the
sum $\sum_{t>1}$ in the equation above needs to be substituted by
$\sum_{t=2,\ldots,\bar{t}}$; for those of second type, this sum must
be substituted by $\sum_{t=\bar{t}+1,\ldots,T}$.

More difficult is the case in which a linear or generalized linear
parameterization is assumed on the transition probabilities. In
order to update the parameters of these models, i.e. the vector
$\b\de$ in (\ref{eq:linear_model}) or (\ref{eq:glm_transition}), we
have to update the corresponding component of the expected value of
the complete data log-likelihood:
\[
\hat{\ell}^*_3(\b\de)=\sum_{t>1}\sum_u\sum_v
\hat{a}^{(t)}_{uv}\log[\pi_{v|u}^{(t)}].
\]
Even in this case, it is natural to apply iterative algorithms, such
as the Fisher-scoring.

For instance, when the multinomial logit parametrization based on
logits of type (\ref{eq:logit_transition}) is assumed, we have
\begin{eqnarray*}
\hat{\b s}_3^*(\b\de)&=&\sum_{t>1}\sum_u [\b V_u^{(t)}]\tr(\hat{\b a}^{(t)}_{u}-n\b\pi_u^{(t)}),\\
\hat{\b F}_3^*(\b\de)&=&n\sum_{t>1}\sum_u [\b
V_u^{(t)}]\tr\{\diag(\b\pi_u^{(t)})-\b\pi_u^{(t)}[\b\pi_u^{(t)}]\tr\}\b
V_u^{(t)},
\end{eqnarray*}
where $\b\pi_u^{(t)}$ is the column vector with $k$ elements  equal
to $\pi_{v|u}^{(t)}$, $v=1,\ldots,k$, and $\hat{\b a}^{(t)}_{u}$ is
the vector of the same dimension with elements $\hat{a}^{(t)}_{uv}$.
In a similar way we can obtain the score vector and the expected
information matrix for other cases, such as that based on global
logits. For  the case of a linear model on the transition
probabilities, or parametrizations of the type illustrated in
Example 2, we refer to \cite{bart:06}.

\subsection{Likelihood ratio testing}\label{sec:testing_constraints}
As usual, in order to test a hypothesis $H_0$ expressed through the
constraints expressed above, we can use the likelihood ratio
statistic
\[
D=-2[\ell(\hat{\b\th})-\ell(\hat{\b\th}_0)],
\]
where $\ell(\hat{\b\th})$ is the maximum value of the likelihood of
the unconstrained model and $\ell(\hat{\b\th}_0)$ is that of the
restricted model.

When the usual regularity conditions hold, the null asymptotic
distribution of the above test statistic is of chi-squared type with
a number of degrees of freedom equal to the number of non-redundant
constraints used to formulate $H_0$. The latter is equal to the
difference in the number of non-redundant parameters between the two
models that are compared. These regularity conditions hold for most
of the constraints formulated in this section. For instance, in the
case of binary response variables, through $D$ we can test the
hypothesis that the conditional probabilities of success follow a
Rasch model; see equation (\ref{eq:rasch}). In this case, $D$ has a
chi-squared null asymptotic distribution with $kT-[T+(k-1)]$ degrees
of freedom, where $T+(k-1)$ is the number of non-redundant
parameters involved in (\ref{eq:rasch}).

The main case where the usual regularity conditions do not hold is
when $H_0$ is formulated by a linear model on the transition
probabilities, see (\ref{eq:linear_model}), such that certain of
these probabilities are constrained to be equal to zero. In this
case, a boundary problem occurs \citep{self:lian:87} and, as proved
by \cite{bart:06}, the null asymptotic distribution is of
chi-bar-squared type, i.e. a mixture of chi-squared distributions
with weights which may be computed by explicit formulae or estimated
by a simple Monte Carlo method; see \cite{shap:88} and
\cite{silv:sen:04}, Chapter 3. An interesting result is when we
assume that the transition matrices depend on only one parameter as
in (\ref{eq:trans_only_one}). In this case, the hypothesis
$H_0:\de=0$ may be tested by the likelihood ratio $D$, whose null
asymptotic distribution corresponding to a mixture between 0 and a
chi-squared distribution with one degree of freedom. The weights of
this mixture are 0.5 and 0.5.
\section{Including individual covariates and relaxing local independence}
\label{sec:cov}
Individual covariates may be included in the measurement model or in
the latent model. Note that the problem of relaxing local
independence is strongly related to that of the inclusion of
individual covariates in the measurement model, because both
extensions rely on suitable parameterizations of the conditional
distribution of the response variables given the latent process.
Once a suitable parametrization has been adopted, local independence
is relaxed by including the lagged response variables, or a suitable
transformation of these variables, among the individual covariates.

Before illustrating the technical details, we need to clarify the
notation. We denote by $\b x^{(t)}_{i}$ the vector of possibly
time-varying individual covariates for subject $i$ at occasion $t$,
$i=1,\ldots,n$, $t=1,\ldots,T$. As mentioned above, this vector may
include the lagged response variables, i.e. $y^{(t-1)}_{i}$ in the
univariate case or $\b y^{(t-1)}_{i}$ in the multivariate case, with
even higher-order lags. The notation for the conditional response
probabilities and initial and transition probabilities must be
changed to take into account the covariates and the lagged
responses. Then, for every $i$, we let
$\phi_{i,y|u}^{(t)}=p(Y^{(t)}_{i}=y|\b x^{(t)}_{i},y^{(t-1)}_{i})$,
$t=1,\ldots,T$, $u=1,\ldots,k$, $y=0,\ldots,l-1$; moreover, we let
$\pi_{i,u}=p(U^{(1)}_{i}=u|\b x^{(1)}_{i})$, $u=1,\ldots,k$, and
$\pi_{i,v|u}^{(t)}=p(U^{(t)}_{i}=v|U^{(t-1)}_{i}=u,\b
x^{(1)}_{i},y^{(t-1)}_{i})$, $t=2,\ldots,T$, $u,v=1,\ldots,k$. In
the multivariate case, we have the probabilities
$\phi_{ij,y|u}^{(t)}=p(Y^{(t)}_{ij}=y|\b x^{(t)}_{i},\b
y^{(t-1)}_{i})$ which are defined, together with the probabilities
$\pi_{i,v|u}^{(t)}$, on the basis of the lagged vector of response
variables $\b y^{(t-1)}_{i}$. Finally, the manifest probability of
the response configuration $\b y$ provided by subject $i$ is denoted
by $f_i(\b y)$ and corresponds to $p(\b Y_i=\b y|\b X_i)$, where $\b
X_i$ is the matrix collecting the vector of covariates $\b
x^{(t)}_{i}$, $t=1,\ldots,T$, for subject $i$.

Before describing in detail the parameterization of the
probabilities on which the measurement and the latent models are
based, an important point to note is that we can still compute in an
efficient way the manifest probabilities $f_i(\b y)$ by using
recursion (\ref{eq:rec}).
\subsection{Covariates in the measurement
model}\label{sec:measurement_cov}
\subsubsection{Univariate case.}
Individual covariates may be included in the measurement model
through a parametrization which recalls that used in (\ref{eq:glm})
to formulate constraints on the conditional distribution of the
response variables given the latent process.

For every subject $i$, let $\b\eta_{i,u}^{(t)}=\b
g(\phi^{(t)}_{i,0|u},\ldots,\phi^{(t)}_{i,l-1|u})$, with $\b
g(\cdot)$ being a link function of the type mentioned in Section
\ref{sec:modelling1}. Then we assume
\begin{equation}
\b\eta_{i,u}^{(t)}=\b Z_{i,u}^{(t)}\b\be,\label{eq:glmcov}
\end{equation}
where $\b Z_{i,u}^{(t)}$ is a design matrix depending on the
covariates in $\b X_i$. We recall that these covariates may include
the lagged response variables when we want to relax local
independence.
\subsubsection{Multivariate case}
In the multivariate case, a parametrization of type
(\ref{eq:glmcov}) can be adopted for each response variable, having
in this way an expression that recalls that in (\ref{eq:glm_multi})
and based, in this case, on design matrices denoted by $\b
Z_{ij,u}$.

An LM model for multivariate data in which the conditional
distribution of the response variables depends on the individual
covariates was formulated by \cite{bart:farc:09}. In this
formulation, the assumption that the response variables in $\b
Y_i^{(t)}$ are conditionally independent given $U_i^{(t)}$ is
relaxed by assuming a marginal parameterization of the type
described in \cite{bart:colo:forc:07}. This implies formulating a
link function with respect to the column vector $\b p_{i,u}^{(t)}$
having elements $p(\b Y_i^{(t)}=\b y^{(t)}|U_i^{(t)}=u,\b
x_i^{(t)})$ for all possible configurations of $\b y^{(t)}$. In
particular, the link functions adopted by \cite{bart:farc:09} may be
formulated as
\begin{equation}
\b\eta^{(t)}_{i,u} =\b C\log[\b M\b p_{i,u}^{(t)}],\label{reparam}
\end{equation}
where $\b C$ and $\b M$ are matrices of simple construction. In
particular, the vector $\b\eta^{(t)}_{i,u}$ contains marginal logits
and marginal log-odds ratios of different types (e.g. local, global,
continuation) for the distribution $\b p_{i,u}^{(t)}$. In this case,
we have a unique design matrix $\b Z_{i,u}^{(t)}$, having then an
expression similar to (\ref{eq:glmcov}), instead of a transition
matrix $\b Z_{ij,u}^{(t)}$ for each response variable. The following
example clarifies the type of parametrization.

\begin{example} {\bf - Marginal model for two response variables.}
Consider the case of $r=2$ variables with two and three levels
($l_1=2$, $l_2=3$, $l_3=3$), which are treated with logits of type
local and global, respectively. Overall, there are 3 logits and 2
log-odds ratios. The logits may be parametrized as follows
\begin{eqnarray*}
\log\frac{\phi_{i1,1|u}^{(t)}}{\phi_{i1,0|u}^{(t)}}&= &\xi_{1u}+[\b x_i^{(t)}]\tr\b\be_1\\
\log\frac{\phi_{i2,1|u}^{(t)}+\phi_{i2,2|u}^{(t)}}{\phi_{i2,0|u}^{(t)}}&= &\xi_{2u}+[\b x_i^{(t)}]\tr\b\be_2\\
\log\frac{\phi_{i2,1|u}^{(t)}}{\phi_{i2,0|u}^{(t)}+\phi_{i2,1|u}^{(t)}}&= &\xi_{3u}+[\b x_i^{(t)}]\tr\b\be_3,\\
\end{eqnarray*}
whereas for the log-odds ratios we have
\[
\log\frac{p(Y_{i1}^{(t)}=1,Y_{i2}\geq y|U_i^{(t)}=u,\b
x_i^{(t)})p(Y_{i1}^{(t)}=0,Y_{i2}<y|U_i^{(t)}=u,\b
x_i^{(t)})}{p(Y_{i1}^{(t)}=0,Y_{i2}\geq y|U_i^{(t)}=u,\b
x_i^{(t)})p(Y_{i1}^{(t)}=1,Y_{i2}\geq y|U_i^{(t)}=u,\b x_i^{(t)})}=
\be_{2+y},
\]
with $y = 1,2$. In this case, the parameter vector $\b\be$ is made
of the subvectors $\b\be_1$, $\b\be_2$ and $\b\be_3$, and the
elements $\be_4$ and $\be_5$.
\end{example}

\subsection{Covariates in the latent model}
\label{sec:covlat}
A natural way to allow the initial and transition probabilities of
the latent Markov chain to depend on the individual covariates is by
adopting a parametrization which recalls that in
(\ref{eq:glm_transition}).

For the initial probabilities, in particular, we assume
\[
\b\la_i = \b V_i\b\de_1,
\]
where $\b\la_i=\b g(\pi_{i,1},\ldots,\pi_{i,k|u})$, $\b V_i$ is a
design matrix depending on the covariates in $\b x_i^{(1)}$ and
$\b\de_1$ is the corresponding parameter vector. Similarly, for the
transition probabilities we have
\begin{equation}
\b\la_{i,u}^{(t)}=\b
V_{i,u}^{(t)}\b\de_2,\label{eq:logit_transition_cov}
\end{equation}
with $\b\la_{i,u}^{(t)}=\b g_u(\pi_{i,1|u},\ldots,\pi_{i,k|u})$,
where the design matrix $\b V_{i,u}^{(t)}$ depends on the covariates
in $\b x_i^{(t)}$ and $\b\de_2$ is the corresponding vector of
parameters. In the above expressions, $\b g(\cdot)$ and $\b
g_u(\cdot)$, $u=1\,\ldots,k$, are link functions that may be
formulated on the basis on different types of logit. Typically, we
use multinomial logits or global logits. As reference category, the
multinomial logits have the first category when modeling the initial
probabilities and category $u$ when modeling the transition
probabilities. In the latter case, an expression similar to
(\ref{eq:logit_transition}) results. Global logits, based on an
expression similar to (\ref{eq:global_transition}) are used when the
latent states are ordered on the basis of a suitable parametrization
of the conditional distribution of the response variables given the
latent process.

Finally, note that we can combine a parametrization of type
(\ref{eq:logit_transition_cov}) with the constraint that certain
transition probabilities are equal to 0 as in
(\ref{eq:triang_transition}).
\subsection{Interpretation of the resulting
models}\label{sec:interpretation}
We introduced two different schemes for including individual
covariates in the model. In formulating an LM model, we suggest to
adopt only one scheme, i.e. to include the covariates in the
measurement model (under the constraints $\pi_{i,u}=\pi_u$,
$\pi_{i,v|u}^{(t)}=\pi_{v|u}^{(t)}$, $i=1,\ldots,n$) or in the
latent model (under the constraint
$\phi_{i,y|u}^{(t)}=\phi_{y|u}^{(t)}$ or
$\phi_{ij,y|u}^{(t)}=\phi_{j,y|u}^{(t)}$, $i=1,\ldots,n$). We advice
against allowing the covriates to affect both the distribution of
the latent process and the conditional distribution of the response
variables given this process. In fact, the two extensions have a
distinct interpretation. Moreover, the resulting LM model is in
general of difficult interpretation and its estimation via EM
algorithm is often cumbersome.

About the model interpretation, we have to clarify that when
covariates are included in the measurement model, the latent
variables are seen as way to account for the unobserved
heterogeneity, i.e. the heterogeneity between subjects that we
cannot explain on the basis of the observable covariates. The
advantage with respect to a standard random effect or latent class
model with covariates is that we admit that the effect of
unobservable covariates has its own dynamics; for further
discussions see \cite{bart:farc:09}.

When the covariates are included in the latent model, we typically
suppose that observeble outcomes indirectly measure a latent trait,
such as the health condition of elderly people, which may evolve
over time. In such a case, the main interest is in modeling the
effect of covariates on the latent trait distribution; see
\cite{bart:lupp:mont:09}.
\subsection{Maximum likelihood estimation}
In the presence of individual covariates, it is convenient to
express the log-likelihood through a function similar to
(\ref{eq:lk}), that is
\[
\ell(\b\th) = \sum_i\log[f_i(\b y_i)],
\]
where $f_i(\b y_i)$ is the manifest probability of the response
configuration provided by subject $i$ given the covariates in $\b
X_i$.

The likelihood function can be maximized by an EM algorithm having a
structure very similar to that outlined in Section \ref{sec:EM}.
This algorithm is based on the complete data log-likelihood that, in
the univariate case, has expression
\begin{eqnarray}
\ell^*(\b\th) &=&
\sum_i\bigg\{\sum_t\sum_u\sum_y\tilde{a}_{i,uy}^{(t)}\log[\phi^{(t)}_{i,y|u}]+\label{eq:comp_lk_cov}\\
&+&\sum_u a_{i,u}^{(1)}\log(\pi_{i,u})+\sum_{t>1}\sum_u\sum_v
a_{i,uv}^{(t)}\log[\pi_{i,v|u}^{(t)}]\bigg\},\nonumber
\end{eqnarray}
where $a_{i,u}^{(t)}$ is a dummy variable equal to 1 if subject $i$
is in latent state $u$ at occasion $t$; with reference to the same
occasion and the same subject, $a_{i,uv}^{(t)}$ is a dummy variable
equal to 1 if this subject moves from state $u$ to state $v$,
whereas the dummy variable $\tilde{a}_{i,uy}^{(t)}$ is equal to 1 if
this subject is in state $u$ and provides response $y$. In the
multivariate case, we have a similar expression for this function,
which depends on the probabilities $\phi^{(t)}_{ij,y|u}$.

The E-step of the algorithm consists of computing the conditional
expected value of the above dummy variables given the observed data
and the current value of the parameters. In practice, we compute
\begin{eqnarray*}
\hat{a}^{(t)}_{i,u}&=&r^{(1)}_i(u|\b y_i),\\
\hat{a}^{(t)}_{i,uv}&=&r_i^{(t)}(u,v|\b y_i),\\
\hat{\tilde{a}}^{(t)}_{i,uy}&=&I(y_i^{(t)}=y)r^{(t)}(u|\b y_i),
\end{eqnarray*}
where $r_i^{(t)}(u|\b y_i)=p(U^{(t)}_{i}=u|\b Y_i=\b y_i)$ and
$r_i^{(t)}(u,v|\b y_i)=p(U^{(t-1)}_i=u,U^{(t)}_i=v|\b Y_i=\b y_i)$;
the conditioning on the observable covariates in $\b X_i$ is
implicit in these expressions. The posterior probabilities may be
computed by applying, subject by subject, the same recursion
illustrated in Section \ref{sec:EM}; see in particular
(\ref{eq:rec2}) and (\ref{eq:defR}).

The M-step consists of maximizing the complete data log-likelihood
expressed as in (\ref{eq:comp_lk_cov}), with each dummy variable
substituted by the corresponding expected value. How to maximize
this function depends on the specific formulation of the model and,
in particular, on whether the covariates are included in the
measurement model or in the latent model. We illustrate this point
below, closely following the scheme adopted in Section
(\ref{sec:EM_cov}).
\subsubsection{Covariates in the measurement model.}
In this case, the initial and transition probabilities of the latent
process are typically assumed to be equal across subjects (see
Section \ref{sec:interpretation}). Then, the common parameters
$\pi_u$ and $\pi_{v|u}^{(t)}$ are updated as
\begin{eqnarray*}
\pi_u &=& \frac{\sum_i\hat{a}_{i,u}^{(1)}}{n},\quad u=1,\ldots,k,\\
\pi_{v|u}^{(t)} &=&
\frac{\sum_i\hat{a}_{i,uv}^{(t)}}{\sum_i\hat{a}_{i,u}^{(t-1)}},\quad
t=2,\ldots,T,\:u,v=1,\ldots,k.
\end{eqnarray*}

On the other hand, under one of the parameterizations illustrated in
Section \ref{sec:measurement_cov}, we have to update the parameter
vector $\b\be$ by an iterative algorithm which maximizes
\[
\hat{\ell}^*_1(\b\be)=
\sum_i\sum_t\sum_u\sum_y\hat{\tilde{a}}_{i,uy}^{(t)}\log[\phi^{(t)}_{i,y|u}]
\]
in the univariate case and a similar expression in the multivariate
case. At this aim we suggest to use an iterative algorithm of
Fisher-scoring type, the same suggested in Section
\ref{sec:EM_measurement_cov}. With binary response variables, in
particular, for the score vector and information matrix we have the
following expressions
\begin{eqnarray*}
\hat{\b s}_1^*(\b\be)&=&\sum_i\sum_t\sum_u\sum_y \hat{\tilde{a}}^{(t)}_{i,uy}(y-\phi^{(t)}_{i,1|u})\b z_{i,u}^{(t)},\\
\hat{\b F}_1^*(\b\be)&=&\sum_i\sum_t\sum_u
\hat{a}^{(t)}_{i,u}\phi^{(t)}_{i,0|u}\phi^{(t)}_{i,1|u}\b
z_{i,u}^{(t)}[\b z_{i,u}^{(t)}]\tr.
\end{eqnarray*}
\subsubsection{Covariates in the latent model.}
When the covariates are included in the latent model, the
conditional distribution of the response variables given the latent
process is the same for all subjects, and depends on the parameters
$\phi_{y|u}^{(t)}$ in the univariate case and on the parameters
$\phi_{j,y|u}^{(t)}$ in the multivariate case. In the first case,
these parameters are updated as
\[
\phi_{y|u}^{(t)}=\frac{\hat{\tilde{a}}_{i,uy}^{(t)}}{\hat{a}_{i,u}^{(t)}},\quad
t=1,\ldots,T,\:u=1,\ldots,k,\:y=0,\ldots,l-1.
\]
A similar rule is adopted in the multivariate case.

Iterative algorithms are instead required to update the parameters
$\b\de_1$ and $\b\de_2$ involved in the latent model, when this
model is formulated as in Section \ref{sec:covlat}. These parameters
are updated by maximizing the functions
\begin{eqnarray*}
\hat{\ell}^*_2(\b\de_1)&=&\sum_i\sum_u
a_{i,u}^{(1)}\log(\pi_{i,u})\\
\hat{\ell}^*_3(\b\de_2)&=&\sum_{t>1}\sum_u\sum_v
a_{i,uv}^{(t)}\log[\pi_{i,v|u}^{(t)}]
\end{eqnarray*}
For the maximization of both functions we again suggest to use the
Fisher-scoring algorithm, that may be implemented along the same
lines as in Section \ref{sec:EM_latent_cov}. The score vectors and
the information matrices exploited by this algorithm may be derived
on the basis of standard rules.
\section{Multilevel extension}\label{sec:multilevel}
We now consider the extension to multilevel data in which subjects
are collected into a given number of clusters and one or more random
effects are used to model the influence of each cluster on the
responses provided by the subjects who are included. The typical
example is that of students collected in classes belonging to
different schools. This extension is strongly related to that at
basis of the mixed LM model
\citep{vand:lang:90,lang:vand:94,kapl:08} and to the extension for
random effects proposed by \cite{Altman:07}.

Even in this case, we may choose to include the random effects in
the measurement model or in the latent model. In the following, we
illustrate the model based on second choice, which we consider more
interesting, closely following what proposed by
\cite{bart:penn:vitt:10}.

In illustrating the multilevel extension, we denote by $H$ the
number of clusters. Every subject in the sample is identified by the
pair of indices $hi$, with $h=1,\ldots,H$ and $i=1,\ldots,n_h$, and
where $n_h$ is the dimension of cluster $h$. Accordingly, we denote
the response of this subject at occasion $t$ by $Y_{hi}^{(t)}$, the
vector of all responses provided by this subject by $\b Y_{hi}$, and
the collection of all responses provided by the subjects in cluster
$h$ by $\b Y_h$. Note that in the multivariate case we have a vector
of responses $\b Y_{hi}^{(t)}=(Y_{hi1}^{(t)},\ldots,Y_{hir}^{(t)})$,
instead of the single variable $Y_{hi}^{(t)}$, where each variable
$Y_{hij}^{(t)}$ has $l_j$ levels, from 0 to $l_j-1$.
\subsection{Model assumptions}
\begin{figure}[h]\centering
\includegraphics[width=11cm]{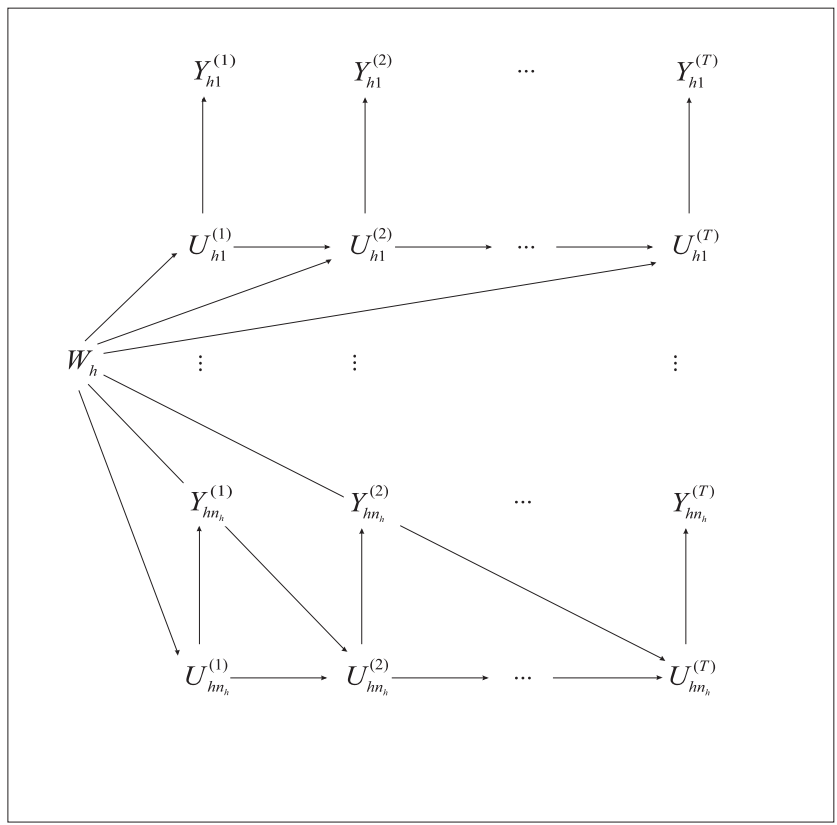}
\caption{\it Path diagram for the multilevel LM
model.}\label{pathmultilevel}
\end{figure}

As usual, the LM model relies on a latent Markov chain $\b
U_{hi}=(U_{hi}^{(1)},\ldots,U_{hi}^{(T)})$ for every subject $hi$.
However, the assumption that these chains are independent, as in the
basic model, is relaxed by assuming that, for every cluster $h$, $\b
U_{h1},\ldots,\b U_{hn_h}$ are conditionally independent given the
latent variable $W_h$. The latter has the role of capturing the
heterogeneity between clusters and, for every $h$, is assumed to
have a discrete distribution with support $\{1,\ldots,m\}$. The
assumption of local independence is retained and formulated by
requiring that each response variable $Y_{hi}^{(t)}$ is
conditionally independent of any other variable in the model
(including the response variables associated to any other subject
and the latent variables at cluster level), given the corresponding
latent state $\b U_{hi}^{(t)}$. All these assumptions are
represented by the path diagram in Figure \ref{pathmultilevel}.

\vspace*{0.25cm}

The above assumptions imply that, for every $h$, the response
vectors $\b Y_{h1},\ldots,\b Y_{hn_h}$ are not independent, but they
are conditionally independent given the latent variable $W_h$.
Marginal independence holds between the vectors of response
variables $\b Y_1,\ldots,\b Y_H$ associated to the different
clusters.

Model parameters are still the conditional response probabilities
and the initial and transition probabilities of the latent process.
The parameters of the first type are still denoted by
$\phi_{y|u}^{(t)}$, $t=1,\ldots,T$, $u=1,\ldots,k$,
$y=0,\ldots,l-1$, and correspond to $p(Y_{hi}^{(t)}=y|U_i^{(t)}=u)$;
a similar notation is adopted in the multivariate case. An extended
notation, which takes into account the multilevel structure, is
instead necessary for the other parameters. In particular, for
$w=1,\ldots,m$, we have $\pi_{u|w} = p(U_{hi}^{(1)}=u|W_h=w)$,
$u=1,\ldots,k$, and $\pi_{v|uw}^{(t)} =
p(U^{(t)}_{i}=v|U^{(t-1)}_{i}=u,W_h=w)$, $t=2,\ldots,T$,
$u,v=1,\ldots,k$, $w=1,\ldots,m$. These parameters correspond,
respectively, to the initial probabilities and the transition
probabilities of the Markov chain associated to a subject belonging
to a cluster with latent variable at level $w$. Finally, we also
have to include the parameters for the distribution of each variable
$W_h$, and then we let $\rho_w=p(W_h=w)$, $w=1,\ldots,m$.

Under the model assumptions, we can express, for each cluster $h$,
the manifest distribution of the vector $\b Y_h$ as
\begin{equation}
f_h(\b y)=p(\b Y_h=\b y)=\sum_w p(\b Y_h=\b
y|W_h=w)\rho_w,\label{eq:manifest_multilevel}
\end{equation}
where $\b y$ has now the same dimension as $\b Y_h$ and is made of
the subvectors $\b y_i$ of the same dimension as $\b Y_{hi}$,
$i=1,\ldots,n_h$. Moreover, we have that
\[
p(\b Y_h=\b y|W_h=w)=\prod_ip(\b Y_{hi}=\b y_i|W_h=w),
\]
where
\begin{eqnarray*}
p(\b Y_{hi}=\b y_i|W_h=w)
&=&\sum_{u^{(1)}}\phi_{y_i^{(1)}|u^{(1)}}\pi_{u^{(1)}}
\sum_{u^{(2)}}\phi_{y_i^{(2)}|u^{(2)}}^{(2)}\pi_{u^{(2)}|u^{(1)}}^{(2)}\cdots\\
&&\cdots\sum_{u^{(T)}}\phi_{y_i^{(T)}|u^{(T)}}^{(T)}\pi_{u^{(T)}|u^{(T-1)}}^{(T)}.
\end{eqnarray*}
This probability may be efficiently computed through recursion
(\ref{eq:rec}).
\subsection{Extension to include cluster and individual covariates and to the multivariate case}
\label{covinittrans}
In the context of multilevel data, covariates are usually available
at cluster and individual levels. Let $\b z_h$ denote the column
vector of available covariates for cluster $h$, let $\b
x_{hi}^{(t)}$ be the column vector of available covariates for
subject $hi$ at occasion $t$, $h=1,\ldots,H$, $i=1,\ldots,n_h$,
$t=1,\ldots,T$. These covariates are included in the model in a way
similar to that illustrated in Section \ref{sec:covlat}. Then we
extend the notation by letting, for each subject $hi$, $\pi_{hi,u|w}
= p(U_{hi}^{(1)}=u|W_h=w,\b x_{hi}^{(1)})$, $u=1,\ldots,k$, and
$\pi_{hi,v|uw}^{(t)} = p(U_{i}^{(t)}=v|U_{i}^{(t-1)}=u,W_h=w,\b
x_{hi}^{(t)})$, $t=2,\ldots,T$, $u,v=1,\ldots,k$, $w=1,\ldots,m$; we
also let $\rho_{h,w}=p(W_h=w|\b z_h)$, $w=1,\ldots,m$, for each
cluster $h$.

Covariates at cluster level are included into the model by adopting
a suitable parametrization for the probabilities $\rho_{h,w}$. A
natural choice is the multinomial logit parametrization, and then we
assume
\begin{equation}
\log\frac{\rho_{h,w}}{\rho_{h,1}}=\ga_{0w}+\b z_h\tr\b\ga_{1w},\quad
w=2,\ldots,m,\label{eq:logit_u}
\end{equation}
where $\b\ga_{12},\ldots,\b\ga_{1m}$ are vectors of regression
coefficients of the same dimension as $\b z_h$, and
$\ga_{02},\ldots,\ga_{0m}$ denote the corresponding intercepts.

A multinomial logit parameterization may also be adopted for the
initial and transition probabilities of the Markov chain associated
to every subject $hi$. However, if the conditional distribution of
the response variables given the latent process is formulated so
that the latent states are ordered, a considerable reduction of
parameters is achieved by adopting a parametrizion based on global
logits. This approach is exploited by \cite{bart:penn:vitt:10}, who
relied on a Rasch parametrization for the distribution of each
response variable given the corresponding latent state. Then, for
the initial probabilities of every Markov chain they assume
\begin{equation}
\log\frac{\pi_{hi,u|w}+\cdots+\pi_{hi,k|w}}{\pi_{hi,1|w}+\cdots+\pi_{hi,u-1|w}}=\de_{0w}+\de_{1u}+(\b
x_{hi}^{(1)})\tr\b\de_2,\quad u=1,\ldots,k,\:
w=2,\ldots,m,\label{eq:logit_ini}
\end{equation}
where $\b\de_2$ is a vector of regression parameters of the same
dimension as $\b x_{hi}^{(t)}$ which is common to every level $u$.
Moreover, the intercepts $\de_{0w}$ depend on the level of $W_h$
(with $\de_{01}\equiv 0$), and the intercepts $\de_{1u}$ depend on
the level of $U_{hi}^{(1)}$. We impose
$\de_{12}\leq\cdots\leq\de_{1m}$ in order to ensure the
invertibility of the global logit parametrization.

Finally, for what concerns the transition probabilities we have
\begin{equation}
\log\frac{\pi_{hi,v+1|uw}^{(t)}+\cdots+\pi_{hi,k|uw}^{(t)}}
{\pi_{hi,1|uw}^{(t)}+\cdots+\pi_{hi,v-1|uw}^{(t)}}=\eta_{0w}+\eta_{1uv}+(\b
x_{hi}^{(t)})\tr\b\eta_2,\label{eq:logit_trans}
\end{equation}
with $u=1,\ldots,k$, $v=1,\ldots,k_2-1$, $w=1,\ldots,m$, and
$t=2,\ldots,T$. As above, $\b\eta_2$ is a vector of regression
coefficients for the individual covariates, the intercepts
$\eta_{0w}$ depend on the level of $W_h$ (with $\eta_{01}\equiv0$),
and the intercepts $\eta_{0uv}$ depend on the levels of
$U_{hi}^{(t-1)}$ and $U_{hi}^{(t)}$ and must be decreasing ordered
in $v$. These parameters may also be allowed to depend on the time
occasion; for instance, we can have $\eta_{0v}^{(t)}$ in place of
$\eta_{0v}$.

Finally note the extension of the model presented to the
multivariate case, where we observe for every $hi$ and $t$ a vector
of responses $\b Y_{hi}^{(t)}$, is straightforward. The basic
assumption to be added is that the variables in this vector are
conditionally independent given the corresponding latent state
$U_{hi}^{(t)}$. \cite{bart:penn:vitt:10} also addressed this case,
allowing for different sets of responses observed at each occasion.
This assumption of conditional independence may also be relaxed by
adopting a parametrization of type (\ref{reparam}), including the
lagged response variables as already discussed in Section
\ref{sec:cov}.
\subsection{Maximum likelihood estimation}
Given an observed sample, the log-likelihood of the model
illustrated above is given by
\[
\ell(\b\th)=\sum_h\log[f_h(\b y_h)],
\]
where $\b y_h$ is the vector containing the responses observed for
all subjects in cluster $h$ and the manifest probability $f(\b y_h)$
is defined in (\ref{eq:manifest_multilevel}).

The EM algorithm is again the mail tool we suggest to maximize
$\ell(\b\th)$. However, its implementation is more difficult with
respect to the versions already illustrated for the presence of
latent variables at two levels. First of all, for the case of
univariate responses and considering the presence of covariates, the
complete data log-likelihood may be expressed as
\begin{eqnarray}
\ell^*(\b\th) &=& \sum_h \sum_w b_{h,w}\log(\rho_{h,w})+\label{eq:complk_multilevel}\\
&+&\sum_h \sum_i\bigg\{ \sum_u a_{hi,u}^{(1)}\log[\pi_{hi,u|w}^{(1)}]+\nonumber\\
&+&\sum_u\sum_v\sum_{t>1}
a_{hi,uv}^{(t)}\log[\pi_{hi,v|uw}^{(t)}]+\nonumber\\
&+&\sum_u\sum_t\sum_y
\tilde{a}_{hi,uy}^{(t)}\log[\phi_{y|u}^{(t)}]\bigg\},\nonumber
\end{eqnarray}
where the dummy variables are defined as in Section \ref{sec:EM_cov}
and $b_{h,w}$ is a dummy variable equal to 1 if cluster $h$ belongs
to latent class $w$ (i.e. $W_h=w$).

At the E-step, we compute the conditional expected value of the
above dummy variables given the observed data and the current value
of the parameters; at the M-step, we maximize the conditional
expected value of $\ell^*(\b\th)$ obtained by substituting each
dummy variable in (\ref{eq:complk_multilevel}) with the
corresponding expected value obtained from the E-step. For details
on the implementation of these steps we refer to
\cite{bart:penn:vitt:10}.

\section{Standard errors, model selection, and path prediction}
\label{sec:se}
In this section we show how to obtain standard errors for the model
parameters, dealing in particular with the information matrix. We
also outline the problem of model selection for what mainly concerns
the choice of the number of latent states. Finally, we outline the
problem of path prediction, i.e. how to find the maximum a
posteriori sequence of latent states for a given subject.
\subsection{Standard errors}
It is well known that the EM algorithm, differently from other
algorithms such as the Newton-Raphson or the Fisher-scoring, does
not directly produce the information matrix of the model. This
matrix is typically used for computing standard errors.

Several methods have been proposed to overcome this difficulty. Most
of these methods have been developed within the literature on hidden
Markov models; for a concise review see \cite{Lyst:Hugh:exac:2002}.
The more interesting methods are based on the information matrix
obtained from the EM algorithm by the technique of
\cite{Loui:find:1982} or related techniques; see for instance
\cite{Turn:Came:Thom:hidd:1998} and \cite{bart:farc:09}. Other
interesting methods obtain the information matrix on the basis of
the second derivative of the manifest probability of the response
variables by a recursion similar to (\ref{eq:rec}); see
\cite{Lyst:Hugh:exac:2002} and \cite{bart:06}.

Among the method related to the EM algorithm, that proposed by
\cite{bart:farc:09} is very simple to implement and requires a small
extra code over that required for the maximum likelihood estimation.
The method exploits a well-known result according to which the score
of the complete data log-likelihood computed at the E-step of this
algorithm corresponds to the score of the incomplete data
log-likelihood. More precisely, we have
\[
\b
s(\b\th)=\frac{\pa\ell(\b\th)}{\pa\b\th}=\left.\frac{\pa\E_{\bar{\bl\th}}[\ell^*(\b\th)|\b
Y]}{\pa\b\th}\right|_{\bar{\bl\th}=\bl\th},
\]
where $\E_{\bar{\bl\th}}[\ell^*(\b\th)|\b Y]$ denotes the
conditional expected value of the complete-data log-likelihood
computed at the parameter value $\bar{\b\th}$. Then, the observed
information matrix, denoted by $\b J(\b\th)$, is obtained as minus
the numerical derivative of $\b s(\b\th)$. The standard error of
each parameter estimate is then obtained as the square root of the
corresponding diagonal element of $\b J(\hat{\b\th})^{-1}$. These
standard errors may be used for hypothesis testing and for obtaining
confidence intervals in the usual way.

The information matrix can also be used for checking local
identifiability at $\hat{\b\th}$. We consider the model to be local
identifiable if the matrix $\b J(\hat{\b\th})$ is of full rank. See
also \cite{mchu:56} and \cite{good:74b}.
\subsection{Model Choice}
In applying a LM model, a fundamental problem is that of the choice
of the number of latent states, denoted by $k$. For multilevel
versions of this model, it is also necessary to choose the number of
support points for the latent variables at cluster level, denoted by
$m$.

In order to choose the above quantities, it is natural to use
information criteria such as the Akaike Information Criterion (AIC),
see \cite{aka:73}, and the Bayesian Information Criterion (BIC), see
\cite{sch:78}. According to first criterion, we choose the number of
states corresponding to the minimum of $AIC =
-2\ell(\hat{\b\th})+2g$, where $g$ is the number of non-redundant
parameters; according to the second, we choose the model with the
smallest value of $BIC = -2\ell(\hat{\b\th})+g\log(n)$.

The performance of the two approaches above have been deeply studied
in the literature on mixture models; see \cite{mcla:peel:00},
Chapter 6. These criteria have also been studied in the hidden
Markov literature for time series, where the two indices above are
penalized with a term depending on the number of time occasions; see
\cite{Boucheron:2007}. From these studies, it emerges that BIC is
usually preferable to AIC, as the latter tends to overestimate the
number of latent states.

The theoretical properties of AIC and BIC applied to the LM models
are less studied. However, BIC is a commonly accepted model choice
criterion even for these models. It has been applied by many
authors, such as \cite{lang:94}, \cite{lang:vand:94}, and
\cite{magi:verm:01}. In particular, \cite{bart:lupp:mont:09}
suggested the use of this criterion together with that of diagnostic
statistics measuring the goodness-of-fit and
goodness-of-classification, whereas a simulation study may be found
in \cite{bart:farc:09}.
\subsection{Path prediction}\label{classmember}
Once the model has been estimated, a relevant issue is that of path
prediction, i.e. finding the most likely sequence of latent states
for a given subject on the basis of the responses he/she provided.
For the $i$th subject in the sample, this is the sequence $\hat{\b
u}_i=(\hat{u}_i^{(1)},\ldots,\hat{u}_i^{(T)})$ such that
\[
p(\hat{\b u}|\b y_i)=\max_{\bl u}p(\b u|\b y_i),
\]
where we use the notation of Section \ref{sec:lm}. This definition
may be simply adapted to more general cases involving individual
covariates and/or multilevel data.

The problem above is different from that of finding the most likely
state occupied by a subject at certain occasion. This problem may be
simply solved on the basis of the posterior probabilities of type
$r^{(t)}(u|\b y)$, whose computation is required within the EM
algorithm; see Section \ref{sec:EM}. In any case, we can use an
efficient algorithm which avoid the evaluation of the posterior
probability $p(\b u|\b y)$ for every configuration $\b u$ of the
latent process. This is known as Viterbi algorithm
\citep{vite:67,juan:rabi:91}, and is illustrated in the following.

For a given subject $i$ with response configuration $\b
y_i=(y_i^{(1)},\ldots,y_i^{(T)})$, let ${\tilde r}_i^{(1)}(u) =
p(U_i^{(1)}=u,y_i^{(1)})$ and, for $t=2,\ldots,T$, let
\[
{\tilde r}_i^{(t)}(u) =\max_{u^{(1)},\ldots,u^{(t-1)}}
p(U_i^{(1)}=u^{(1)},\ldots,U_i^{(t-1)}=u^{(t-1)},U_i^{(t)}=u,
y_i^{(1)},\ldots,y_i^{(t)}).
\]
A forward recursion can be used to compute the above quantities, and
a backward recursion based on these quantities can then be used for
path prediction:
\begin{enumerate}
\item for $u=1,\ldots,k$ compute $\tilde{r}^{(1)}_{i}(i)$ as $\pi_up(y^{(1)}_{i}|U^{(1)}_{i}=u)$;
\item for $t=2,\ldots,T$ and $v=1,\ldots,k$ compute
${\tilde r}^{(t)}_{i}(v)$ as
\[
p(y_i^{(t)}|U^{(t)}_{i}=v)\:\max_u [{\tilde
r}^{(t-1)}_{i}(u)\pi_{uv}^{(t)}];
\]
\item find the optimal state $\tilde{u}^{(T)}_{i}$
as $\hat{u}^{(T)}_{i}=\arg\max_u \tilde{r}^{(T)}_{i}(u)$;
\item for $t=T-1,\ldots,1$,
find $\hat{u}_i^{(t)}$ as $\tilde{u}^{(t)}_{i}=\arg\max_u \tilde
{r}^{(t)}_{i}(u)\pi_{u\hat{u}_i^{(t+1)}}$.
\end{enumerate}

All the above quantities are computed on the basis of the ML
estimate of the parameter $\b\th$ of the model of interest.

\section{Empirical illustrations}
\label{applications}
In order to illustrate the approaches reviewed in this paper, we
provide a synthetic overview of some interesting applications
appeared in the literature. We describe different datasets, the
aspects involved in practically fitting and using the LM models, and
the interpretation of the results.
\subsection{Marijuana Consumption dataset}
The univariate version of the LM model described in Section
\ref{sec:uni} was applied by \cite{bart:06} to analyze a marijuana
consumption dataset based on five annual waves of the ``National
Youth Survey" \citep{ell:huiz:mena:89}. The dataset concerns $n=237$
respondents who were aged 13 years in 1976.  The use of marijuana
was measured through five ordinal variables, one for each annual
wave, with three categories corresponding to: ``never in the past",
``no more than once in a month in the past year", and ``once a month
in the past year". The substantive research question is whether
there is an increase of marijuana use with age.

Using BIC, \cite{bart:06} selected an LM model with three latent
states, homogeneous transition probabilities, and a parsimonious
parameterization for the measurement model based on global logits,
which recalls that in (\ref{eq:rasch}). This parametrization is
based on one parameter for each latent state, which may be
interpreted as the tendency to use marijuana for a subject in this
state, and one cutpoint for each response category. Then, the latent
states may be ordered representing subjects with ``no tendency to
use marijuana", ``incidental users of marijuana", and ``subjects
with high tendency to use marijuana". Also note that the cutpoints
are common to all the response variables, since these variables
correspond to repeated measurements of the same phenomenon under the
same circumstances. This is because the dynamics of the marijuana
consumption is only ascribed to the evolution of the underlying
tendency of this consumption.

\cite{bart:06} also tested different hypotheses on the transition
matrix of the latent process. In particular, he found that the
hypothesis that the transition matrix has a tridiagonal structure,
i.e.
\[
\b\Pi^{(t)}=\pmatrix{1-\de_1 & \de_1 & 0\cr
                     \de_2 & 1-(\de_2+\de_3) & \de_3\cr
                     0 & 1-\de_4 & \de_4},
\]
cannot be rejected. This hypothesis implies that the transition from
state $u$ to latent state $v$ is only possible when $v=u-1$ or $v=
u+1$.

The results under the selected LM model say that, at the beginning
of the period of observation, the 89.6\% of the sample is in the
first class (lowest tendency to marijuana consumption) and the 1.5\%
is in the third class (highest tendency to marijuana consumption).
An interesting interpretation of the pattern of consumption emerges
from the estimated transition matrix. A large percentage of subjects
remains in the same latent class, but almost 25\% of accidental
users switches to the class of high frequency users. From the
estimated marginal probabilities of the latent classes emerge that
the tendency to use marijuana increases with age, since the
probability of the third class increases across time.
\subsection{Educational dataset}
An interesting illustration of an LM model with individual
covariates was given by \cite{verm:lang:bock:99}. They used data
from an educational panel study conducted by the ``Institute for
Science Education in Kiel (Germany)" \citep{hoff:lehr:todt:85}. A
cohort of secondary school pupils was interviewed once a year from
grade 7 to grade 9 with respect to their interests in physics as a
school subject. The response variables have been dichotomized with
categories ``low" and ``high" to avoid sparseness of the observed
frequency table. Based on these data, the LM model is used to draw
conclusions on whether interest in physics depends on the interest
in the previous period of observation and on two available
covariates: sex and grade in physics at the present time.

\cite{verm:lang:bock:99} estimated a univariate LM model with both
initial and transition probabilities of the latent process depending
on the available covariates according to a multinomial logit
parametrization. In this model, the measurement error was
constrained to be the same for all time points, meaning that the
conditional distribution of the response variables given the latent
state is the same for every occasion. Then, they relaxed some of the
basic assumptions of the LM model, such as the assumption that the
Markov chain is first-order.

According to the parameter estimates of the selected model, there is
a significant effect of sex and grade on the interest in physics.
Pupils with higher grades are more interested in physics than pupils
with lower grades, girls are less interested in physics than boys.
Moreover, the interest has a positive effect on the grade at the
next time occasion. For the boys, the probability of switching from
``low" to ``high interest" is larger than that for girls, as well as
to keep their interest high.
\subsection{Criminal dataset}
The multivariate version of the LM model where both the initial and
the transition probabilities of the latent process depend on
time-constant covariates was illustrated by
\cite{bart:penn:fran:07}. They analyzed the conviction histories of
a cohort of offenders who were born in England and Wales in 1953.
The offenders were followed from the age of criminal responsibility,
10 years, until the end of 1993. They were grouped in 10 major
categories and gender was included in the model as explanatory
variable. The analysis was based on $T=6$ age bands: 10-15, 16-20,
21-25, 26-30, 31-35 and 36-40 years.

The adopted LM model allows to estimate trajectories for
behavioral types which are determined by the criminal conviction
grouping. It also allows to give rise to a general population
sample by augmenting the observed sample with not-convicted
subjects.

According to \cite{bart:penn:fran:07}, the fit of the model is
considerably improved by relaxing the assumption of homogeneity
of the latent Markov chain, but retaining the constraint that males
and females have the same transition probabilities. In particular,
they selected a model based on partially homogeneity, as described
in Section \ref{sec:modelling2}; equation
(\ref{eq:partial_homo}). Therefore there are two transition
probability matrices: the first for transitions up to time $\bar{t}$
and the second from time $\bar{t}$ and beyond. The choice of
$\bar{t} = 2$ has been made on the basis of the BIC.

In summary, the selected model is based on a partially homogeneous
Markov chain with five latent states, different initial and equal
transition probabilities for males and females. From the estimated
conditional probabilities of conviction for any offence group and
any latent state, it was possible to determine classes of criminal
activity. In accordance to the typologies found in the
criminological latent class literature, these classes are
interpreted as: ``non-offenders", ``incidental offenders", ``violent
offenders", ``theft and fraud offenders" and ``high frequency and
varied offenders".

The estimated initial probabilities show that in the first age band
the percentage of males who are incidental offenders is higher than
that of females. The common estimated transition probabilities for
males and females from age band 10-15 to age band 16-20, and from
one age band to the others  for offenders over 16, show that the
first transition occurs at an early age, 16 years, which in western
society represents the peak of the age-crime curve.

At the first time occasions, ``incidental offenders" have a quite
high probability of persistence when moving from the age band 10-15
to age band 16-20. Moreover, ``theft and fraud offenders" are mainly
females and they have a high chance of moving to the class of
non-offenders. The ``high frequency and varied offenders" are mainly
males and they have a high persistence.

The estimated transition probabilities from age band 16-20 to the
subsequent age bands show that the subjects belonging to the latent
state of ``non-offenders" have a very low chance of becoming
offenders; ``theft and fraud offenders" and ``violent offenders"
have a high probability of dropping out of crime. From the estimated
proportion of males and females in each latent state at every time
occasion it can be seen that 7\% of males are ``violent offenders"
at age 16-20 years and 32\% are ``incidental offenders" at the same
age. Only 3\% of females are ``theft and fraud offenders" at age
16-20 years.
\subsection{Dataset on financial products preferences}
An interesting analysis of data obtained from face to face
interviews of the household ownership of 12 financial products was
offered by  \cite{pass:verm:bijm:07}.  The panel was conducted by a
market research company among 7676 Dutch households in 1996, 1998,
2000 and 2002. To have an accurate representation of the products
portfolio, the households were asked to retrieve their bank and
insurance records in order to check which product they owned.
Households that dropped out were replaced to ensure the
representativeness of the sample for the population with respect to
demographic variables, such as age, income and marital status.

The aim of the study was to get insights on the developments of the
individual household product portfolio and the effect of demographic
covariates on such development. It also concentrates on predicting
future behaviors of acquisition.

The authors proposed to use a time homogeneous multivariate LM
model, with time-varyng covariates affecting the latent process as
in Section \ref{sec:covlat}. They added additional assumptions to
the model, such as constant conditional probabilities of the
response variables given the latent process. This is done to avoid
manifest changes, so that the product penetration levels are
consistent in latent states over measurement occasions. Moreover,
they formulated the model with a time-constant effect of the
covariates on the transition probabilities.

The model selected on the basis of BIC is based on nine latent
states. These states can be ordered according to increasing
penetration levels across the analyzed products, which range from
bonds, the most commonly owned product, to saving accounts. The
results highlight some divergences from common order of
acquisitions, such as the acquisition of a mortgage before owing a
credit card or viceversa. Loans and unemployment insurance are
most often acquired. According to the estimated transition
matrix there is a high persistence in the same latent state: only
14\% of the households changed latent state in the period of the
study. The most common switch is from latent class 7 to 8, where
latent state 7 is characterized by the acquisition of mortgage, life
insurance, pension fund, car insurance, and saving
accounts, whereas latent state 8  for all the previous products plus
the credit card. Another common switch is from latent state 4 to 7,
where the first is characterized by the acquisition of life
insurance, pension fund, car insurance, and savings accounts. This
means that multiple products were acquired between consecutive
measurement occasions.

Income, age of the head of the household, and household size have a
significant effect on the initial and transition probabilities
according to the Wald test. The covariate values implying a larger
probability of belonging to an initial latent class also imply a
greater probability of switching into the same latent class. For
example larger households are relatively often found in latent
states where overall product ownership probabilities are relatively
low.

The prediction of future purchase of a financial product was
performed on the basis the posterior latent state membership
probabilities for each household at the last occasion $T$, given all
other observed information. To assess the accuracy of the
forecasting, the authors used the Gini coefficient as a measure of
concentration. Considering the empirical results in the last wave
referred to year 2002, which was not considered when estimating the
model, the authors showed that, for most products, the prediction
equations are effective for forecasting household acquisition.
\subsection{Job position dataset}
The multivariate LM model with covariates affecting the manifest
probabilities proposed by \cite{bart:farc:09} was applied by these
authors to data extracted from the ``Panel Study of Income Dynamics"
database (University of Michigan). These data concern $n=1446$ women
who were followed from 1987 to 1993. The binary response variables
are fertility, indicating whether a woman had given birth to a child
in a certain year and  employment, indicating whether she was
employed. The covariates are:  race (dummy variable equal to 1 for a
black woman),  age (in 1986),  education (year of schooling), child
1-2 (number of children in the family aged between 1 and 2 years,
referred to the previous year),  child 3-5,  child 6-13, child 14-,
income of the husband (in dollars, referred to the previous year).

The main issue concerns the direct effect of fertility on
employment. Also of interest are the strength of the state
dependence effect for both response variables and how these
variables depend on the covariates.

\cite{bart:farc:09} used an LM model with covariates affecting the
manifest probabilities since they were interested in separately
estimating the effect of each covariate on each outcome. The
proposed LM model allows to separate these effects from the
unobserved heterogeneity, by modeling the latter with a latent
Markov process. In this way, unobserved heterogeneity effects on the
response variables are allowed to be time-varying; this is not
allowed neither within a LC model with covariates nor in the most
common random effect models.

The model selected using AIC and BIC is with three latent states.
Under this model, race has a significant effect on fertility, but
not on employment according to the estimates of the parameters
affecting the marginal logits of fertility and employment and the
log-odds ratio between these variables. Age has a stronger effect on
fertility than on employment. Education has a significant effect on
both fertility and employment, whereas the number of children in the
family strongly affects only the first response variable and income
of the husband strongly affects only the second one.

The log-odds ratio between the two response variables, given the
latent state, is negative and highly significant, meaning that the
response variables are negatively associated when referred to the
same year. On the other hand, lagged fertility has a significant
negative effect on both response variables and lagged employment has
a significant effect, which is positive, on both response variables.
Therefore, fertility has a negative effect on the probability of
having a job position in the same year of the birth, and the
following one. Employment is serially positively correlated (as
consequence of the state dependence effect) and fertility is
negatively serially correlated.

From the estimates of the support points for each latent state it
may be deduced that the latent states correspond to different levels
of propensity to give birth to a child and to have a job position.
The first latent state corresponds to subjects with the highest
propensity to fertility and the lowest propensity to have a job
position. On the contrary, the third latent state corresponds to
subjects with the lowest propensity to fertility and the highest
propensity to have a job position. Finally, the second state is
associated to intermediate levels of both propensities. The two
propensities are negatively correlated.

Overall, it results that the 78.5\% of women started and persisted
in the same latent state for the entire period, whereas for the
21.5\% of women had one or more transitions between states. The
presence of these transitions is in accordance to the rejection of
the hypothesis that a LC model is suitable for these data.
\subsection{Dataset on anorectic patients}
An interesting extension of the model to account for a hierarchical
structure has been  recently proposed by \cite{rijm:vans:boec:07}.
They illustrated the model by a  novel application using a data set
from an ecological momentary assessment study
\citep{vans:rijm:piet:vand:07} on the course of emotions among
anorectic patients. At nine occasions for each of the seven days of
observations, 32 females with eating disorders received a signal and
were asked to rate themselves on a 7-point scale with respect to the
intensity with which they experienced 12 emotional states. These
were taken from the following emotional categories: ``anger and
irritation", ``shame and guilt", ``anxiety and tension", ``sadness
and loneliness", ``happiness and joy", ``love and appreciation". The
response has been dichotomized (0-2 vs. 3-6) and the signal has been
considered equally spaced. The aim of the study was to detect the
course of emotion among the patients.

As a preliminary analysis, \cite{rijm:vans:boec:07} used the
univariate version of the LM model without covariates. They treated
each person by day combination as a separate case assuming that the
data stemming from different days were independent and that the
parameters were constant over days. On the basis of such model the
authors selected four latent states. The first state is interpreted
as positive mood, the third as negative mood, the second as low
intensity for all emotions except tension, and the fourth as neutral
to moderately positive mood. According to the estimated transition
matrix there is high persistence in the same state. The probability
of moving from state 1 to state 2 is 0.14, from state 3 to 4 is 0.18
and from state 4 to 3 is 0.14 . They noted that there is an indirect
transition from state 3 to state 1 via the emotionally more neutral
state 4. Over the days there is an increase of the marginal
probabilities of states 1 and 4 indicating that the mood of patients
tends to become better later on in the day.

\cite{rijm:vans:boec:07} also used a hierarchical LM model by
introducing a latent variable at day level to account for the fact
that data stemming from different days are not independent. They
modeled the transition between latent states at day and signal
levels by a first-order time homogeneous Markov chain. They
estimated a model with two states at the day level and two signal
states within each day-state. For the first day state, the signal
state 1 is characterized by high probabilities of experiencing
positive emotions and low probabilities of negative emotions.
Therefore this state is interpreted as positive mood. Instead the
second day state is interpreted as negative mood. In the signal
state 2 positive emotions are not well separated from negative
emotions and the state is considered as an emotionally neutral to
moderately positive state. A tendency to experience more positive
emotions emerges from the estimated initial and conditional
probabilities of the chain over day and signal.
\subsection{Dataset on student math achievement}
A multilevel version of the latent Markov model was applied by
\cite{bart:penn:vitt:10} to analyze how the cognitive math
achievement changes over time of schooling. Multilevel models are
standard tools for the analysis of such kind of data as they allow
to explicitely consider a hierarchical structure such as students
nested in classes and schools.

The data derived from the repeated administration of test items to
students attending public and non-public middle schools in an
Italian Region. A set of dichotomously scored items was administered
at the end of each of the three years of school (28 at the end of
the first year, 30 at the end of the second, and 39 at the end of
the third year) to 1,246 students who progressed from Grade 6 to
Grade 8. They are from 13 public and 7 non-public middle schools.
The recorded social background characteristics of the students are
father's and mother's education. The school characteristics are:
number of students, number of teachers, students-teachers ratio, and
years since school opened.

The main issue concerns the evaluation of the performance of
schools, measured in terms of achievement attained by pupils at the
end of the period of formal schooling. The proposed approach takes
into account imperfect persistence of achievement, heterogeneity in
learning and measurement error in test scores contributing to get an
unbiased estimate of the added value of public and non-public
schools.

On these data the authors fitted a multilevel LM model for
multivariate data with a structure similar to that presented in
Section \ref{sec:multilevel}. On the basis of the BIC, they selected
the model based on $m=4$ support points for the latent variable at
cluster level and $k=6$ states for the latent Markov chain at
individual level. Each latent state is associated to an estimated
the ability level, so that these states can easily identify the
group of the most proficient students and that of the least
proficient students. The ordered estimated conditional probabilities
for each latent class and each set of items administered at each
grade show that the difficulty of the items is increasing over time
and they allow to identify the items which are tailored to
distinguish less capable students.

The estimates of the intercepts and the regression coefficients for
the logistic model at cluster level based on parametrization
(\ref{eq:logit_u}) allowed interpretation of the four cluster latent
classes. The class of cluster $D$ is mainly characterized by the
non-public schools with less than 18 years since school opened and
with a ratio between students and teachers smaller than 8. As the
value of the covariate ``ratio between student and teachers"
increases there is more chance that the class is of type A rather
than of type B. As the value of the year of activity of the school
increases there is more chance that the class is of type A rather
than of type B.

The estimated parameters of the global logits defined on the initial
and transition probabilities are interpreted on the basis of
formulae (\ref{eq:logit_ini}) and (\ref{eq:logit_trans}). According
to them the classes of those schools belonging to type $B$ are less
helpful in increasing math ability in the first year of the middle
school, as compared to the classes of the other clusters. The same
classes contribute less also from Grade 6 to Grade 7 but they
contribute the most to student's math ability from Grade 7 to 8. The
ability of the students increases for those having highly educated
fathers. The magnitude of this increase is stronger at Grade 6
compared to the other grades. They identify how the covariates
related to each class affect the probability that this class is of
type A, B, C or D. The most interesting conclusion is that classes
of students in public schools are mostly of type A (32\%), B (23\%),
and C (38\%), whereas those in non-public schools are only of type A
(79\%) and D (21\%). The estimated regression coefficients related
to the individual covariates show that the ability of the students
increases for those having higher educated fathers.
\section{Conclusions and further developments}
In this paper we presented a review of the latent Markov model
which starts with the basic model of \cite{wigg:73}, for both
univariate and multivariate data. We also illustrated extensions to
the inclusion of covariates and to multilevel data, and outlined
maximum likelihood estimation and related inferential methods. We
tried to keep the presentation simple and provide references to
methods which are of easy implementation, taking ideas from methods
widely applied in the literature on hidden Markov models for
time series. Moreover, in order to retain a simple presentation
we dealt with the case of balanced data in which all subjects are
observed at the same number of occasions, the same response
variables are obtained at every occasion, and there are no missing
responses.

The potentialities of the LM approach are illustrated by summarizing
the results of several applications available in the literature. It
is worth noting that the number of applications available in the
literature is growing; recent applications also involve situations
which are not typically those of longitudinal data. Examples concern
the estimation of closed-populations on the basis of
capture-recapture data \citep{bart:penn:07} and the assessment of
the peer review process of grant applications \citep{Bornmann:08}.

In this section we would like to summarize further developments of
the LM models that, in our opinion, deserve particular attention:
\begin{itemize}
\item {\em More flexible parametrizations for the conditional distribution
of the response variables given the latent process}: We refer in
particular to parameterizations which address the problem of
multidimensionality when the model is applied, for instance, to data
coming from the administration of test items assessing different
types of ability. This is a problem well known in the Item Response
Theory, which may concern different fields of applications such as
psychology and assessment of the quality of life. Flexible
parameterizations are also of interest when the response variables
are of mixed nature; for instance, we may have a mixture of binary,
ordinal, and continous variables.
\item {\em Higher-order Markov chains}: the assumptions that the latent
Markov chain if of first-order may be restrictive when there is
``memory effect" such that the latent state of a subject at a given
occasion depends on two or more previous occasions. The extension to
second or higher-order Markov chains
is rather simple to implement and has been already applied
in certain contexts. However, to our knowledge, a systematic
illustration of this higher-order approach does not seem to be
available in the literature.
\item {\em Missing responses}: It is well known that longitudinal
datasets typically suffer from the problem of missing responses due
to several reasons. In applying an LM model, a trivial solution to
this problem consists in eliminating from the sample subjects with
at least one missing response. Obviously, this solution may lead to
a strong bias of the parameter estimates when the assumption that
these responses are missing at random is not plausible. The interest
here is in formulating LM models in which the event of missing
response is explicitly modeled given the latent state and the
evolution of the latent process takes into account this possibility
by adopting a suitable parametrization of the transition matrices.
\item {\em Multilevel data with evolution at cluster level}: We
presented in Section \ref{sec:multilevel} a multilevel extension of
the LM model in which only one latent variable is used to represent
the effect of each cluster. However, this effect may have its own
dynamics as the effect at individual level. Hence, of interest is
the extension in which a sequence of occasion-specific latent
variables is used to model the effect of each cluster. The resulting
model may then be formulated by assuming the existence of two nested
latent Markov chains. Given the chain at cluster level, another
chain is associated to each individual in the cluster. This model if
of particular interest, but it seems difficult to provide a simple
method for maximum likelihood estimation, given the difficulty of
obtaining the manifest distribution of the response variables.
\end{itemize}

\section*{Acknowledgements}
The authors are supported by the EIEF research grant "Advances in
non-linear panel models with socio-economic applications". F.
Bartolucci and F. Pennoni also acknowledge the financial support
from PRIN 2007.

\bibliographystyle{apalike}
\bibliography{reference}
\end{document}